\documentclass[11pt,a4paper,leqno]{amsart}

\usepackage{amsmath,amssymb,amsfonts,amsthm,euscript} 

\def\dis
{\displaystyle}
\def\eps
{\varepsilon}
\def\e
{\varepsilon}
\def\si{{\sigma}}
\def\om{{\omega}}

\def\C{{\mathbb C}}
\def\R{{\mathbb R}}
\def\N{{\mathbb N}}
\def\Z{{\mathbb Z}}

\def\Q{{\mathbb Q}}

\def\Tend#1#2{\mathop{\longrightarrow}\limits_{#1\rightarrow#2}}

\def\d{{\partial}}

\theoremstyle{plain}
\newtheorem{theo}{Theorem}[section]
\newtheorem{lem}[theo]{Lemma}
\newtheorem{cor}[theo]{Corollary}
\newtheorem{prop}[theo]{Proposition}
\newtheorem{hyp}[theo]{Assumption}

\theoremstyle{definition}
\newtheorem{defin}[theo]{{\sc Definition}}

\theoremstyle{remark}
\newtheorem*{rema}{{\sc Remark}}

\numberwithin{equation}{section}

\begin{document}
\title[Semiclassical NLS with
potential]{Semiclassical Nonlinear Schr\"odinger 
equations with potential and focusing initial data} 
\author[R. Carles]{R{\'e}mi Carles}
\address[R. Carles]{MAB, UMR CNRS 5466\\
Universit{\'e} Bordeaux 1\\ 351 cours de la Lib{\'e}ration\\ 33~405 Talence
cedex\\ France}
\email{carles@math.u-bordeaux.fr}
\author[L. Miller]{Luc Miller}
\address[L. Miller]{Centre de Math{\'e}matiques, UMR CNRS 7640\\ 
{\'E}cole Polytechnique\\ 91~128
Palaiseau cedex\\ France}
\address{and: \'Equipe MODAL'X\\ UFR SEGMI \\ B\^atiment G\\ 200
Avenue de la R\'epublique\\ 92~001 Nanterre cedex
\\ France}
\email{miller@math.polytechnique.fr}
\thanks{This work was partially supported by the ACI grant
  ``{\'E}quation des ondes : oscillations, dispersion et
contr{\^o}le''. These results were improved while the first author was
in the University of Osaka, invited by N.~Hayashi, to whom he wishes
to express his gratitude.}
%\begin{abstract}
%We study the asymptotic behaviour of solutions to semi-classical
%nonlinear Schr\"odinger equations with a 
%potential, for concentrating and oscillating initial data, when the
%nonlinearity is 
%repulsive and the potential is a polynomial of degree at most two. We
%describe the separate roles of the nonlinearity and of 
%the potential, with tools which seem to be specific to the class of
%potentials that we consider. We also discuss the case of more general
%subquadratic potentials. 
%\end{abstract}
\maketitle

\section{Introduction}
\label{sec:intro}

We study the semi-classical limit $\e \to 0$ of
solutions ${\bf u}^\e:(t,y)\in \R\times \R^n\to \C$ of the equation
\begin{equation*}
i\e \d_t {\bf u}^\e +\frac{1}{2}\e^2 \Delta {\bf u}^\e = V(y){\bf u}^\e
+ \lambda |{\bf u}^\e|^{2\si} {\bf u}^\e\ ,
\end{equation*}
where $\lambda >0$ (the nonlinearity is repulsive), with concentrating
initial data 
\begin{equation*}
{\bf u}^\e (0,y)= R\left(
\frac{y-y_0}{\e}\right)e^{i\frac{y\cdot \eta_0}{\e}}\ .
\end{equation*}
Similar problems were studied for attractive nonlinearities 
($\lambda <0$), by Bronski and Jerrard (\cite{BJ00}), and
Keraani (\cite{Keraani02}). In that case, if the power is
$L^2$-subcritical ($\si <2/n$) and $R$ is the ground state
solution of an associated scalar elliptic equation, then when $V$ is
smooth with $V \in W^{2,\infty}$, the following asymptotics holds in
$X:=L^\infty_{loc}(\R;L^2(\R^n))$, 
\begin{equation}\label{eq:bj}
\begin{split}
& \frac{1}{\e^{n/2}}\left\|{\bf u}^\e(t,y) -R\left(
\frac{y-y(t)+\e y^\e(t)}{\e}\right)
e^{i\frac{y\cdot \eta(t)}{\e}+i\theta^\e(t)}
\right\|_X=O\left(\sqrt\e\right) , \\
& \frac{1}{\e^{n/2}}\left\|\e \nabla_y\left({\bf u}^\e (t,y)-R\left(
\frac{y-y(t)+\e y^\e(t)}{\e}\right)
e^{i\frac{y\cdot \eta(t)}{\e}+i\theta^\e(t)}\right)
\right\|_X=O\left(\sqrt\e\right) ,
\end{split}
\end{equation}
where $\theta^\e(t) \in [0,2\pi[$, $y^\e:\R \to \R^n$ is locally
uniformly bounded and $(y(t),\eta(t))$ are the integral
curves associated to the classical Hamiltonian
\begin{equation}\label{eq:hamil}
p(t,y,\tau,\eta) = \tau +\frac{1}{2}|\eta|^2 +V(y),
\end{equation}
with initial data $(y_0,\eta_0)$. 

In this paper, we address the case of a defocusing nonlinearity
($\lambda >0$), when the potential is a polynomial of degree at most
two. 

In the case $\lambda >0$, a
different qualitative behaviour is expected. Intuitively, dispersive
effects prevent the solution from keeping a concentrating aspect as
in \eqref{eq:bj}, for it is well known (see e.g. \cite{Caz}) that the
solutions to the nonlinear Schr\"odinger equation
\begin{equation}\label{eq:nls}
i\d_t \psi +\frac{1}{2}\Delta \psi = |\psi|^{2\si}\psi\ ,
\end{equation}
have the same dispersive properties as the solutions to the linear
Schr\"odinger equation, under suitable assumptions on $\si$ and
$\psi(0,y)$. In the case where the potential $V$ is the harmonic
potential, $V(y) = \om^2 |y|^2$, it was proved in \cite{CaIHP} that
when $y_0=\eta_0=0$,  
the nonlinear term is relevant so long as the dispersive effects are
not too strong. This is so in a boundary layer of size $\e$. Past this
boundary layer, the nonlinear term becomes negligible, and the
potential $V$ imposes the dynamical behaviour of the solution. In the
case of an \emph{isotropic} potential, 
\begin{equation}
\label{eq:harmo}
V(y)=\frac{1}{2}\sum_{j=1}^n \om_j^2 y_j^2, 
\end{equation}
where all the $\om_j$'s are equal, then focusing at the origin occurs
at times $t=k \pi$, $k\in \Z$, and each focus crossing is described in
terms of the Maslov index (this phenomenon is linear) and the nonlinear
scattering operator associated to \eqref{eq:nls}. The case $\eta_0=0$,
$y_0\in \R^n$, is also discussed, and we explain below how to infer
the more general case $(y_0,\eta_0)\in \R^n\times \R^n$ (see \eqref{eq:chgt}). 

The case where the $\om_j$'s are (all positive) not necessarily equal
is also discussed in \cite{CaIHP}. The conclusion is that the nonlinear
term is not relevant outside the initial boundary layer if and only if
two of the $\omega_{j}$'s are rationnaly independent. In the present
paper, we consider the 
case of a generalized quadratic potential which excludes this case.

More precisely, we assume that the potential $V$ is of the form
\begin{equation}\label{def:V1}
V(y)=\sum_{1\leq j,k\leq n}\alpha_{jk}y_j y_k + \sum_{j=1}^n \beta_j
y_j +\gamma,
\end{equation}
where the constants $a_{jk}$, $b_j$ and $c$ are real. We first notice
that up to changing the origin and the basis, we can assume that the
potential has a more rigid form. 
\begin{lem}\label{lem:reduc}
Let $V$ given by \eqref{def:V1}. 
There exist $\widehat y \in \R^n$, and a family $f_1,\ldots,f_n\in
\R^n$ of orthogonal unit vectors such that, with $\widehat y$ as a
new origin, the potential $V$ writes, in 
the basis $(f_1,\ldots,f_n)$, 
\begin{equation*}
V(x)= \frac{1}{2}\sum_{j=1}^n \delta_j \om_j^2 x_j^2 + \sum_{j=1}^n b_j
x_j +c,
\end{equation*}
where $\om_j >0$, $\delta_j \in \{-1,0,1\}$, $b_j$, $c\in \R$ and for
every $j$, $\delta_j b_j=0$. The real numbers 
$$\frac{1}{2}\delta_j \om_j^2\ ,\ \ j=1\ldots n\ ,$$
are the eigenvalues of the quadratic part of $V$. 
\end{lem}
\begin{proof}
Consider the quadratic part of the potential $V$, 
\begin{equation*}
q(y)=\sum_{1\leq j,k\leq n}\alpha_{jk}y_j y_k.
\end{equation*}
It is well-known that there exists a family $f_1,\ldots,f_n\in
\R^n$ of orthogonal unit vectors such that, in this new basis, $q$
writes
\begin{equation*}
q(\widetilde y)=\frac{1}{2}\sum_{j=1}^n \delta_j \om_j^2 \widetilde y_j^2,
\end{equation*}
where $\om_j >0$, $\delta_j \in \{-1,0,1\}$. In this basis, $V$ is of
the form
\begin{equation*}
V(\widetilde y)= \frac{1}{2}\sum_{j=1}^n \delta_j \om_j^2 \widetilde
y_j^2 + \sum_{j=1}^n \widetilde \beta_j 
\widetilde y_j +\gamma,
\end{equation*}
with $\widetilde \beta_j  \in \R$. If $\delta_j=0$, we take $b_j
=\widetilde \beta_j $, and if $ \delta_j\not =0$, we use the
one-dimensional formula,
\begin{equation*}
x^2 +2a x = \left( x+a \right)^2 -a^2. 
\end{equation*}
The lemma follows. 
\end{proof}
In these new coordinates, the Laplace operator is not changed, and 
the initial value problem we are interested
in becomes
\begin{equation}\left\{
\begin{split}
i\e \d_t {\bf u}^\e +\frac{1}{2}\e^2 \Delta {\bf u}^\e &= V(x){\bf u}^\e
+ \lambda |{\bf u}^\e|^{2\si} {\bf u}^\e\ ,\\
{\bf u}^\e (0,x)&= R\left(
\frac{x-x_0}{\e}\right)e^{i\frac{x\cdot \xi_0}{\e}}e^{i\frac{\kappa}{\e}}\ ,
\end{split}\right.
\end{equation}
for some $x_0,\xi_0\in \R^n$, $\kappa \in \R$. Notice that
$\widetilde{\bf u}^\e$, defined by $\widetilde{\bf u}^\e(t,x):=
{\bf u}^\e(t,x)e^{i(ct + \kappa)/\e}$, solves
\begin{equation}\label{eq:intermed}\left\{
\begin{split}
i\e \d_t \widetilde{\bf u}^\e +\frac{1}{2}\e^2 \Delta 
\widetilde{\bf u}^\e &= \left(V(x)-c\right)\widetilde{\bf u}^\e
+ \lambda |\widetilde{\bf u}^\e|^{2\si} \widetilde{\bf u}^\e\ ,\\
\widetilde{\bf u}^\e (0,x)&= R\left(
\frac{x-x_0}{\e}\right)e^{i\frac{x\cdot \xi_0}{\e}}\ ,
\end{split}\right.
\end{equation}
We can thus assume $c=0$. We make an additional assumption on the
potential. 
\begin{hyp}
\label{hyp:pot}
We suppose that the potential satisfies the following properties.\\
1. It is of the form
\begin{equation}
\label{eq:form}
V(x)=\frac{1}{2}\sum_{j=1}^n \delta_j \om_j^2 x_j^2 + \sum_{j=1}^n b_j
x_j,
\end{equation}
where $\om_j >0$, $\delta_j \in \{-1,0,1\}$, $b_j$, $c\in \R$ and for
every $j$, $\delta_j b_j=0$. \\
2. Either there exists $j$ such that $\delta_j \not = 1$, or
$\delta_j=1$ for all $j$ and the $\om_j$'s are not pairwise rationally
dependent: 
$$\exists j\not = k,\ \frac{\om_j}{\om_k}\not \in \Q.$$ 
\end{hyp}
\begin{rema}
We allow negative coefficients for the potential (case
$\delta_j =-1$). In that case, the energy of ${\bf u}^\e$ which is formally
independent of time,
\begin{equation}\label{eq:energy}
E^\e = \frac{1}{2}\|\e \nabla_x {\bf u}^\e(t)\|_{L^2}^2 +\frac{1}{\si
+1}\|{\bf u}^\e (t)\|_{L^{2\si+2}}^{2\si+2} +\int V(x) |{\bf
u}^\e(t,x)|^2dx\ \ , 
\end{equation}
contains negative terms which are not controlled by the positive terms
(in particular, by the $H^1$-norm). Therefore, even the issue of
global existence in $H^1$ is not obvious. We prove that for any $T>0$,
${\bf u}^\e $ cannot blow up for $|t|\leq T$, provided that $\e$ is
sufficiently small ($0<\e \leq \e(T)$). Notice that in the case of an
isotropic negative quadratic potential ($\delta_j =-1$ and $\om_j=\om$
for all $j$), global existence for fixed $\e$ was proved in
\cite{CaRep}. 
\end{rema}

Assumption~\ref{hyp:pot} has a simple geometric consequence. Forget
the nonlinear term for a moment, and consider the classical
Hamiltonian $p$ given by \eqref{eq:hamil}. Because $V$ is of the form
given by \eqref{eq:form}, the bicharacteristic curves starting from
any point $(x_0,\xi_0)\in \R^n\times\R^n$ can be computed explicitly. 
They solve
the differential equation
\begin{equation*}\left\{
\begin{split}
\dot t =1\ ; &\ \dot x(t)=\xi(t)\ , \\
\dot \tau =0\ ; & \ \dot \xi(t)=-\nabla V\left(x(t)\right)\ ,\\
x(0)=x_0\ ; & \ \xi(0)=\xi_0\ .
\end{split}\right.
\end{equation*}
Introduce the auxiliary functions,
\begin{equation}\label{eq:gh}
g_j(t)=\left\{
\begin{split}
\frac{\sin (\om_j t)}{\om_j}\ ,\ &\textrm{ if }\delta_j=1\ ,\\
t\ ,\ &\textrm{ if }\delta_j=0\ ,\\
\frac{\sinh (\om_j t)}{\om_j}\ ,\ &\textrm{ if }\delta_j=-1\ .
\end{split}
\right.
\ \ \ ; \ \ 
h_j(t)=\left\{
\begin{split}
\cos(\om_j t)\ ,\ &\textrm{ if }\delta_j=1\ ,\\
1\ ,\ &\textrm{ if }\delta_j=0\ ,\\
\cosh (\om_j t)\ ,\ &\textrm{ if }\delta_j=-1\ .
\end{split}
\right.
\end{equation}
Then the bicharacteristic curves are
given by
\begin{equation}\label{eq:bicar}
%\begin{split}
x_j(t)= h_j(t)x_{0j} + g_j(t)\xi_{0j}-\frac{1}{2}b_j t^2  \ ;\ \
\xi_j(t)= h_j(t)\xi_{0j} - \delta_j \om_j^2 g_j(t)x_{0j}-b_j t.
%\end{split}
\end{equation}
As the analysis will prove later on, the second part of
Assumption~\ref{hyp:pot} implies that except at time $t=0$, the energy
is never concentrated at one point. Some new concentrations may happen
for $t\not =0$ (if $\delta_j=1$ for at least one $j$), but on a vector
space of dimension at least one, for which the nonlinear term turns
out to be subcritical in the limit $\e\to 0$.  

First, assume that $x_0=\xi_0=0$. Taking $u^\e :=
\e^{-n/2}\lambda^{1/(2\si)}\widetilde{\bf u}^\e$ as a new unknown 
turns \eqref{eq:intermed} into
\begin{equation}\label{eq:pb}
\left\{
\begin{split}
i\e \d_t u^\e +\frac{1}{2}\e^2 \Delta u^\e &= V(x)u^\e
+ \e^{n\si} |u^\e|^{2\si} u^\e\ ,\\
u^\e_{\mid t=0}&= \frac{1}{\e^{n/2}}\varphi\left(
\frac{x}{\e}\right) \ ,
\end{split}
\right.
\end{equation}
where $\varphi$ is given by $\varphi :=\lambda^{1/(2\si)} R$. 
As we mentioned already, we expect the caustic crossing at time $t=0$
to be described by the scattering operator associated to
\eqref{eq:nls}. For this operator to be well-defined, we make a
second assumption, on the initial datum and the nonlinearity.

\begin{hyp}\label{hyp:nl}
The initial datum $\varphi$ and the power $\si$ are such that: \\
1. $\dis \varphi \in \Sigma :=\left\{ f\in H^1(\R^n)\ ; \ |x|f \in
L^2(\R^n) \right\}$, where $\Sigma$ is equipped with the norm
$$\|f\|_\Sigma = \|f\|_{L^2} +\|\nabla f\|_{L^2} +\|xf\|_{L^2} \ .$$ 
2. $1\leq n\leq 5$ and $\si >1/2$, so that the nonlinearity
$|z|^{2\si}z$ is twice differentiable.\\
3. If $n=1$, we assume in addition $\si >1$.\\
4. If $3\leq n\leq 5$, we take $\si <\frac{2}{n-2}$. \\
5. If $n\leq 2$, we assume
\begin{itemize}
\item Either $\si >\dis\frac{2-n+\sqrt{n^2+12n+4}}{4n}$,
\item Or $\|\varphi\|_\Sigma\leq \delta$ sufficiently small. 
\end{itemize}
\end{hyp}
\begin{rema}
The assumption $\varphi \in \Sigma$ makes the energy \eqref{eq:energy}
well defined at time $t=0$.
\end{rema}
\begin{rema}
The assumption $\si <\frac{2}{n-2}$
is needed for a complete $H^1$ theory on \eqref{eq:nls} to be available (see
e.g. \cite{Caz}). The assumption $\si >1/2$, used later on for the
nonlinearity to be twice differentiable, therefore imposes the
restriction $n\leq 5$. 
\end{rema}
\begin{rema}
The third and fifth points of the above assumption are here to insure
the existence of a complete scattering theory for
\eqref{eq:nls}. When $n\geq 3$, this theory is available because $\si
>1/2$. Denote $U_0(t)=e^{i\frac{t}{2}\Delta}$ the free
Schr\"odinger group. From \cite{HT87} and \cite{CW92}, since $\varphi
\in \Sigma$, there exist $\psi_\pm\in\Sigma$ such that the unique
solution $\psi$ to \eqref{eq:nls} such that $\psi_{\mid t=0}=\varphi$
satisfies
\begin{equation}\label{eq:scatt}
\lim_{t\to \pm \infty}\left\| U_0(-t)\psi(t)-\psi_\pm \right\|_\Sigma
=0 \ . 
\end{equation}
\end{rema}
We can now state our main result in the case $x_0=\xi_0=0$. 

\begin{theo}\label{theo}
Suppose that Assumptions~\ref{hyp:pot} and \ref{hyp:nl} are
satisfied. \\
1. For any $T>0$, there exists $\e(T)>0$ such that for $0<\e\leq
   \e(T)$, \eqref{eq:pb} has a unique solution $u^\e \in
   C([-T,T];\Sigma)$. \\
2. This solution satisfies the following asymptotics.
\begin{itemize}
\item For any $\Lambda >0$,
\begin{equation}\label{eq:asym1}
\begin{split}
\limsup_{\e \to 0}\sup_{|t|\leq \Lambda \e}\Big(&
\left\|u^\e(t)-v^\e(t) \right\|_{L^2}+
\left\|\e\nabla_x u^\e(t)-\e\nabla_x v^\e(t) \right\|_{L^2}\\
& +
\left\|\frac{x}{\e}u^\e(t)- \frac{x}{\e}v^\e(t) \right\|_{L^2}\Big)
   =0\ ,
\end{split}
\end{equation}
where 
\begin{equation}\label{eq:v}
v^\e(t,x)=\frac{1}{\e^{n/2}}\psi \left(\frac{t}{\e},\frac{x}{\e} \right), 
\end{equation}
and $\psi\in C(\R;\Sigma)$ is the solution to \eqref{eq:nls} such that
$\psi_{\mid t=0}=\varphi$. 
\item Beyond this boundary layer, we have
\begin{equation}\label{eq:asym2}
\begin{split}
\limsup_{\e \to 0}\sup_{\Lambda \e \leq \pm t\leq T}\Big(&
\left\|u^\e(t)-u_\pm^\e(t) \right\|_{L^2}+
\left\|\e\nabla_x u^\e(t)-\e\nabla_x u_\pm^\e(t) \right\|_{L^2}\\
&+\left\|x u^\e(t)-x u_\pm^\e(t) \right\|_{L^2}\Big)
   \Tend \Lambda {+\infty} 0,
\end{split}
\end{equation}
where $u^\e_\pm \in C(\R;\Sigma)$ are the solutions to 
\begin{equation}\label{eq:upm}\left\{
\begin{split}
i\eps \d_t u_\pm^\e +\frac{1}{2}\e^2\Delta u_\pm ^\eps &=
V(x)u_\pm^\e\ ,  \\
u^\e_{\pm \mid t=0} &= \frac{1}{\e^{n/2}}\psi_\pm
\left( 
\frac{x}{\e}\right),
\end{split}\right.
\end{equation}
and $\psi_\pm$ are given by \eqref{eq:scatt}. 
\end{itemize}
\end{theo}
\begin{rema}
This result can be viewed as a nonlinear analog to a result due to
Nier. In \cite{NierENS} (see also \cite{NierX}), the author studies the
problem 
\begin{equation}\label{eq:nier}
\left\{
\begin{split}
i\eps\d_t u^\eps +\frac{1}{2}\eps^2\Delta u^\eps &
=V(x)u^\eps +U\left( \frac{x}{\eps}\right)u^\eps,\\
u^\eps_{\mid t=0}& = \frac{1}{\eps^{n/2}}\varphi\left( \frac{x}{\eps}\right),
\end{split}
\right.
\end{equation}
where $U$ is a short range potential. The potential $V$ in that case
is bounded as well as
all its derivatives. In that paper, the author proves that under
suitable assumptions, the influence of $U$ occurs near $t=0$ and is
localized near the origin, while only the value $V(0)$ of $V$ at the origin
is relevant in this r\'egime. For times $\eps \ll |t|<T_*$, the
situation is different: the potential $U$ becomes negligible, while
$V$ dictates the propagation. As in our paper, the transition between
these two r\'egimes is measured by the scattering operator associated
to $U$. 

Assumption~\ref{hyp:nl} implies in particular $n \sigma >1$, which 
 makes the nonlinear term short range. With
our scaling for the nonlinearity, this perturbation is relevant only
near the focus, where the potential is negligible, while
the opposite occurs for $\eps\ll |t|\leq T$. 
\end{rema}
%\begin{rema}
%We will see in Section~\ref{sec:general} that the above theorem still
%holds if $V(x)$ is replaced by $V(x)+\e {\tt V}(t,x)$, where ${\tt
%V}(t,x)$ is a smooth potential, subquadratic in $x$. 
%\end{rema}
The case $x_0=\xi_0=0$ turns out not to be so particular in the case
of a potential $V$ satisfying \eqref{eq:form}, when no linear term is
present, that is $b_j=0$, $\forall j$. Introduce the
change of variables
\begin{equation}
\label{eq:chgt}
\begin{split}
{\tt u}^\e(t,x) &= u^\e(t,x-x(t))e^{iS(t,x)/\e}\ ,\\
\textrm{ with }S(t,x)&= x\cdot
\xi(t)-\frac{1}{2}\big(x(t)\cdot \xi(t) -x_0\cdot \xi_0 
\big)\ ,
\end{split}
\end{equation}
where $x(t)$ and $\xi(t)$ are given by \eqref{eq:bicar}. 
It is easy to check that if $u^\e$ solves \eqref{eq:pb} with
$x_0=\xi_0=0$, then ${\tt u}^\e$ solves 
 \begin{equation}\label{eq:pb2}
\left\{
\begin{split}
i\e \d_t {\tt u}^\e +\frac{1}{2}\e^2 \Delta {\tt u}^\e &= V(x){\tt u}^\e
+ \e^{n\si} |{\tt u}^\e|^{2\si} {\tt u}^\e\ ,\\
{\tt u}^\e_{\mid t=0}&= \frac{1}{\e^{n/2}}\varphi\left(
\frac{x-x_0}{\e}\right)e^{ix\cdot \xi_0/\e} \ .
\end{split}
\right.
\end{equation}
\begin{cor}\label{cor:origin}
Let $(x_0,\xi_0)\in \R^n\times \R^n$. 
Under Assumptions~\ref{hyp:pot} and \ref{hyp:nl}, with $b_j=0$,
$\forall j$,
 we have: \\
1. For any $T>0$, there exists $\e(T)>0$ such that for $0<\e\leq
   \e(T)$, \eqref{eq:pb2} has a unique solution ${\tt u}^\e \in
   C([-T,T];\Sigma)$. \\
2. This solution satisfies the following asymptotics.
\begin{itemize}
\item For any $\Lambda >0$,
\begin{equation}\label{eq:asym1'}
\begin{split}
\limsup_{\e \to 0}\sup_{|t|\leq \Lambda \e}\Bigg(&
\left\|{\tt u}^\e(t)-{\tt v}^\e(t) \right\|_{L^2}+
\left\|\e\nabla_x {\tt u}^\e(t)-\e\nabla_x {\tt v}^\e(t)
\right\|_{L^2} \\
& +\left\|\frac{x-x(t)}{\e} {\tt u}^\e(t)- \frac{x-x(t)}{\e}{\tt v}^\e(t)
\right\|_{L^2} \Bigg)
   =0\ ,
\end{split}
\end{equation}
where 
\begin{equation*}
{\tt v}^\e(t,x)=\frac{1}{\e^{n/2}}\psi \left(\frac{t}{\e},\frac{x-x(t)}{\e}
\right)e^{iS(t,x)/\e},  
\end{equation*}
$\psi\in C(\R;\Sigma)$ is the solution to \eqref{eq:nls} such that
$\psi_{\mid t=0}=\varphi$ and $S$ is given by \eqref{eq:chgt}. 
\item Beyond this boundary layer, we have
\begin{equation}\label{eq:asym2'}
\begin{split}
\limsup_{\e \to 0}\sup_{\Lambda \e \leq \pm t\leq T}\Bigg(&
\left\|{\tt u}^\e(t)-{\tt u}_\pm^\e(t) \right\|_{L^2}+
\left\|\e\nabla_x {\tt u}^\e(t)-\e\nabla_x {\tt u}_\pm^\e(t)
\right\|_{L^2} \\
&+\left\|\left(x-x(t)\right)\left({\tt u}^\e(t)-
 {\tt u}_\pm^\e(t)\right)
\right\|_{L^2}  \Bigg) 
   \Tend \Lambda {+\infty} 0,
\end{split}
\end{equation}
where ${\tt u}^\e_\pm \in C(\R;\Sigma)$ are the solutions to 
\begin{equation*}\left\{
\begin{split}
i\eps \d_t {\tt u}_\pm^\e +\frac{1}{2}\e^2\Delta {\tt u}_\pm ^\eps &=
V(x){\tt u}_\pm^\e\ ,  \\
{\tt u}^\e_{\pm \mid t=0} &= \frac{1}{\e^{n/2}}\psi_\pm
\left( 
\frac{x-x_0}{\e}\right)e^{ix\cdot \xi_0/ \e},
\end{split}\right.
\end{equation*}
and $\psi_\pm$ are given by \eqref{eq:scatt}. 
\end{itemize}
\end{cor}
\begin{rema}
The functions ${\tt u}^\e_\pm$ are also given by 
${\tt u}^\e_\pm(t,x) = u^\e_\pm(t,x-x(t))e^{iS(t,x)/\e}$. 
\end{rema}

\begin{rema}
The change of variable \eqref{eq:chgt} could also be used in the case
of an isotropic (attractive) harmonic potential to generalize the
results of \cite{CaIHP}. 
\end{rema}

\begin{rema}
The above corollary shows in particular that the results stated in
Theorem~\ref{theo} are independent of the fact that the concentrating
point is a critical point for the potential $V$. 
\end{rema}

\begin{rema}
After this article was written, it was noticed that we can go further
into reducing the assumptions. Denote $b=(b_1,\ldots,b_n)$, and define
$u^\e_\sharp$ by
\begin{equation*}
u^\e_\sharp(t,x):= u^\e\left(t, x -\frac{t^2}{2}b \right)e^{i\left( 
tb\cdot x -\frac{t^3}{3}|b|^2\right)/\e}\ .
\end{equation*}
As noticed in \cite{CN}, if $u^\e$ solves \eqref{eq:pb}, then
$u^\e_\sharp$ solves the same 
initial value problem, with $V$ replaced by 
$$V_\sharp(x) = V(x) -b\cdot x\ ,$$
which satisfies Assumption~\ref{hyp:pot} and has no linear
part. Therefore, the conclusions of Corollary~\ref{cor:origin} still
hold without the assumption $b=0$. 
\end{rema}

This paper is organized as follows. In Section~\ref{sec:linear}, we
study the linear equations \eqref{eq:upm}. We introduce some tools
which are relevant in the nonlinear setting, and prove that under
Assumption~\ref{hyp:pot}, possible refocusings occur with less
intensity for $t\not =0$ than for $t=0$. In Section~\ref{sec:loc}, we
establish local existence results in $\Sigma$ for \eqref{eq:pb} when
$\e$ is fixed, for general subquadratic potentials. In
Section~\ref{sec:cl}, we prove the first asymptotics of
Theorem~\ref{theo}, and the proof of Theorem~\ref{theo} is completed
in Section~\ref{sec:beyond}. Finally, we examine in
Section~\ref{sec:general} the asymptotic behaviour of $u^\e$ solution
to \eqref{eq:pb} when $V$ is a general subquadratic potential, not
necessarily of the form \eqref{eq:form}.
\section{The linear equation}\label{sec:linear}

In this section, we analyze some properties of solutions of the
equation
\begin{equation}\label{eq:lin}
i\e \d_t u^\e +\frac{1}{2}\e^2 \Delta u^\e = V(x)u^\e\ .
\end{equation}
Under Assumption~\ref{hyp:pot}, it turns out that some tools which are
classical in a linear setting (Heisenberg observables) are very
helpful to study nonlinear problems. Introduce the unitary group
\begin{equation}\label{eq:group}
U^\e(t):= \exp i\frac{t}{\e} \left( \frac{\e^2}{2}\Delta -V(x) \right). 
\end{equation}
This group is well-defined for subquadratic potentials (see
\cite{ReedSimon2}, p.~199), and in particular under our assumptions. 

We consider the following Heisenberg observables (see e.g. \cite{Robert}),
\begin{equation}\label{eq:obs}
A_1^\eps(t):= U^\e(t)\frac{x}{\e}U^\e(-t)\ \ ; \ \ \ A_2^\eps(t):=
U^\e(t)i\e \nabla_x U^\e(-t).
\end{equation}
They solve
\begin{equation*}
%\begin{split}
\d_t A_1^\eps(t)=U^\e(t)i \nabla_x U^\e(-t)=\frac{1}{\e}A_2^\eps(t)\
;\ \ 
\d_t A_2^\eps(t)=-U^\e(t) \nabla_xV U^\e(-t).
%\end{split}
\end{equation*}
Therefore,
\begin{equation*}
\begin{split}
\d_t^2 A_{1,j}^\eps(t)&=-\frac{1}{\e}U^\e(t) \d_j V U^\e(-t)\\
&=-\delta_j \om_j^2 U^\e(t)\frac{x_j}{\e}U^\e(-t)
-\frac{b_j}{\e}\\
&=-\delta_j \om_j^2 A_{1,j}^\eps(t)-\frac{b_j}{\e} . 
\end{split}
\end{equation*}
We thus have explicitly,
\begin{equation}\label{eq:opgene}
\begin{split}
A_{1,j}^\e (t):=& \frac{x_j}{\e}h_j(t)+ ig_j(t)\d_j\ -\frac{b_j}{2\e}t^2,\\
A_{2,j}^\e (t):=& -\delta_j \om_j^2 x_j g_j(t) +ih_j(t)\e \d_j -b_j t .
\end{split}
\end{equation}
These operators inherit interesting properties which we list below. 
\begin{lem}\label{lem:ab}
The operators $A_{\ell,j}^\e$ satisfy the following properties.
\begin{itemize}
\item They commute with the linear part of \eqref{eq:pb},
\begin{equation}\label{eq:commut}
\left[A_{\ell,j}^\e(t), i\e \d_t +\frac{1}{2}\e^2\Delta -V(x)
\right]=0\ ,\ \  \forall (\ell,j)\in \{1,2\}\times \{1,\ldots,n\}\ .
\end{equation}
\item Denote 
\begin{equation*}
\begin{split}
\phi_1(t,x)&:= \frac{1}{2}\sum_{k=1}^n
\left(\frac{h_k(t)}{g_k(t)} x_k^2 -b_k t x_k-\frac{t^3}{12}b_k^2\right)\ ,\\
\phi_2(t,x)&:= -\frac{1}{2}\sum_{k=1}^n\left(
\delta_k \om_k^2 \frac{g_k(t)}{h_k(t)} x_k^2 +2 b_k
t x_k +\frac{t^3}{3}b_k^2\right)\ .
\end{split}
\end{equation*}
Then $\phi_1$ and $\phi_2$ are well-defined for almost every $t$, and
\begin{equation}\label{eq:factor}
\begin{split}
A_{1,j}^\e(t) & =ig_j(t) e^{i\phi_1(t,x)/\e} \d_j
\left(e^{-i\phi_1(t,x)/\e} \ \cdot \right), \\  
A_{2,j}^\e(t) &=i\e h_j (t) e^{i\phi_2(t,x)/\e} \d_j
\left(e^{-i\phi_2(t,x)/\e} \ \cdot \right).
\end{split}
\end{equation}
\item For $r\geq 2$, and $r< \frac{2n}{n-2}$ if $n\geq 3$ ($r\leq \infty$ if
$n=1$), define $\delta(r)$ by
$$\delta(r)\equiv n\left(\frac{1}{2}-\frac{1}{r}\right).$$
Define $P^\e(t)$ by 
\begin{equation*}
P^\e(t):= \prod_{j=1}^n \Big(|g_j(t)| +\e| h_j(t)|\Big)^{1/n}.
\end{equation*}
There exists $C_r$ such
that, for any $f\in \Sigma$,
\begin{equation}\label{eq:sobolev}
\|f\|_{L^r}  \leq \frac{C_r}{P^\e(t)^{\delta(r)}}
\|f\|_{L^2}^{1-\delta(r)}
\max_{\ell,j }\|A_{\ell,j}^\e(t)f\|_{L^2}^{\delta(r)}.
\end{equation}
\item For any function $F\in C^1(\C, \C)$ satisfying the gauge
invariance condition 
$$\exists G\in C(\R_+,\R),\ F(z)=zG(|z|^2),$$
one has, for any $(\ell,j)\in
\{1,2\}\times\{1,\ldots,n\}$ and  almost all $t$, 
\begin{equation}\label{eq:jauge}
A_{\ell,j}^\e(t)F(w)= \d_zF(w)A_{\ell,j}^\e(t)w - \d_{\bar
z}F(w)\overline{A_{\ell,j}^\e(t)w}.
\end{equation} 
\end{itemize}
\end{lem}
\begin{proof}
The first point follows the definition of Heisenberg observables (Von
Neumann equation). 
The second is straightforward computation.
The third point is a
consequence of the well-known Gagliardo-Nirenberg inequalities, and of
\eqref{eq:factor}. The last point is also a consequence of
\eqref{eq:factor}. 
\end{proof}
\begin{rema}
In the definition of $\phi_1$ (resp. $\phi_2$), the factor $t^3
  b_k^2/12 $ (resp.   $t^3 b_k^2/3$) may seem
artificial, for it plays no role in the formula \eqref{eq:factor}. We
introduced these terms because their presence implies that $\phi_1$
  and $\phi_2$
solve the eikonal equation
$$\d_t \phi + \frac{1}{2}|\nabla_x \phi |^2 +V(x)=0\ .$$
This point is discussed further in details in Section~\ref{sec:eik}. 
\end{rema}
\begin{rema}
As noticed in \cite{CaRep}, the fact that our operators enjoy the
properties to be Heisenberg observables \emph{and} factorized as in
\eqref{eq:factor} is due to Assumption~\ref{hyp:pot}. We prove in
Section~\ref{sec:general} that other potentials cannot meet these two
properties.  
\end{rema}

To conclude this section, we explain why the second point of
Assumption~\ref{hyp:pot} implies that there is no ``strong'' focusing
outside $t=0$ for \eqref{eq:pb}. As we will see in the proof of
Theorem~\ref{theo}, this is so because the solutions to \eqref{eq:upm}
do not concentrate at one single point for $t\not = 0$. 

Let $(\ell,j)\in \{1,2\}\times\{1,\ldots,n\}$. Because of
\eqref{eq:commut}, $A^\e_{\ell,j}u_\pm^\e$ solve \eqref{eq:lin}, and
$$\left\|A^\e_{\ell,j}u_\pm^\e (t) \right\|_{L^2} =
\left\|A^\e_{\ell,j}u_\pm^\e (0) \right\|_{L^2}=O(1)\ , \ \textrm{ as
} \e \to 0. $$ 
Thus, for any $r$ as in Lemma~\ref{lem:ab}, there exists $C$
independent of $\e$ and $t$ such that, 
\begin{equation*}
\| u_\pm^\e (t)\|_{L^r} \leq \frac{C}{P^\e(t)^{\delta(r)}}\ \ . 
\end{equation*}
Notice that the concentration of $u_\pm^\e$ is equivalent to the
cancellation of the $g_j$'s. Assume that exactly $p$ functions $g_j$'s
cancel at 
time $t_0$. For the
corresponding $h_j$'s, we have $h_j(t_0)=1$, and $P^\e(t_0)$ is of
order \emph{exactly} $\e^{p/n}$ as $\e$ goes to zero.  The functions
$u_\pm^\e$ concentrate on a space of dimension $n-p$. 

At time $t=0$, we have
\begin{equation}\label{eq:upm0}
\| u_\pm^\e (0)\|^r_{L^r}= \frac{1}{\e^{nr/2}}\int
\left| \psi_\pm \left( \frac{x}{\e}\right)\right|^r dx
=O\left(\e^{-r\delta(r)} \right).
\end{equation}
From the second point of
Assumption~\ref{hyp:pot}, if for $t_0\not =0$, $p$ functions $g_j$'s
cancel, then necessarily, $p<n$, and
\begin{equation}\label{eq:upmt0}
\| u_\pm^\e (t_0)\|_{L^r} =O\left(\e^{-\delta(r)p/n}\right). 
\end{equation}
Comparing \eqref{eq:upm0} and \eqref{eq:upmt0} (recall that $p<n$)
shows that the amplification of the $L^r$-norms cannot be so strong as
at time $t=0$. Since the scaling for the nonlinear term in
\eqref{eq:pb} is critical for the concentration at one point, it is
subcritical for any other concentration, this is why the nonlinear
term is relevant only near the origin in the asymptotics stated in
Theorem~\ref{theo}. This heuristic argument is made rigorous in
Section~\ref{sec:beyond}, and uses the following lemma.
\begin{lem}\label{lem:P}
Let $V$ satisfy Assumption~\ref{hyp:pot}, and denote ${\underline \om}
=\min \om_j$. Let $\delta>0$ and $k>1$ such that $\delta k>1$. 
Then 
\begin{equation*}
\limsup_{\e \to 0}
\e^{-\frac{1}{k}+\delta }\left( \int_{\Lambda \e}^{\pi/(2{\underline \om})}
\frac{dt}{P^\e(t)^{\delta k}}
\right)^{1/k}\Tend \Lambda {+\infty} 0\ .
\end{equation*}
Moreover, for any
$T>0$, there exists $C>0$ independent of $\e\in ]0,1]$, such that
\begin{equation*}
\left( \int_{\pi/(2{\underline \om})}^T
\frac{dt}{P^\e(t)^{\delta k}}
\right)^{1/k}\leq C \e^{\frac{1}{k}-\delta +\frac{\delta}{n}}\ .
\end{equation*}
\end{lem}
\begin{proof}[Sketch of the proof] The functions $g_j$'s may cancel at
times $m\pi/\om_j$, for $m\in \Z$. 
For $t\in [\Lambda \e, \pi/(2{\underline \om})]$, 
$$P^\e(t) \geq \frac{C}{t},$$
and the first part of the lemma follows. For the second part,
split the considered integral into a sum
of the form
\begin{equation*}
\int_{\pi/(2{\underline \om})}^{\pi/{\underline \om} -\e} + 
\int_{\pi/{\underline \om}-\e}^{\pi/{\underline \om} +\e} +
\int_{\pi/{\underline \om} +\e}^{\pi/{\om_j} -\e} + \ldots
+ \int_{\pi/\om_l +\e}^{T}.
\end{equation*}
We noticed that if at time $m\pi/\om_j$, $g_j$ cancels, then at most
$n-1$ functions $g_l$'s cancel, and
$$P^\e(t)\geq C \e^{-1 +1/n}\ ,\ \ \forall t\in \left[
\frac{m\pi}{\om_j}-\e, \frac{m\pi}{\om_j}+\e\right]\ .$$
This shows that integrals of the form 
$$\int_{m\pi/{\om_j}-\e}^{m\pi/{\om_j} +\e}$$
yield the announced estimate. Other integrals are estimates in a
similar fashion. 
\end{proof}

\section{Local existence results}\label{sec:loc}

In this section, we establish local existence results for nonlinear
Schr\"odinger equations with a general subquadratic potential. This is
a natural generalization of \eqref{eq:pb}, and will be needed in
Section~\ref{sec:general}. Consider a potential ${\tt V}$ satisfying
the following properties. 

\begin{hyp}\label{hyp:gen}
The potential ${\tt V}: \R\times \R^n \mapsto \R$ depends on $t$ and $x$, and
satisfies:\\ 
1.  For fixed $t$, ${\tt V}(t,.)\in
C^\infty( \R^n ,\R)$. We also assume that ${\tt V}$ is a measurable function
of $(t,x)\in  \R\times \R^n$. \\
2. For $\alpha \in \N^n$, define
$$M_\alpha (t)=\sup_{x \in \R^n}|\d_x^\alpha {\tt V}(t,x)| + \sup_{x \leq 1
}|{ \tt V}(t,x)| .$$
We assume that for any multi-index satisfying $|\alpha|\geq 2$,
$M_\alpha \in L^\infty_{loc}(\R)$. 
\end{hyp}
Notice that the first point of Assumption~\ref{hyp:pot} implies
Assumption~\ref{hyp:gen}. Denote
\begin{equation}\label{def:ttU}
{\tt U}^\e(t) := \exp\left(i\frac{t}{\e}\left(\frac{\e^2}{2}\Delta
-{\tt V}\right) \right).
\end{equation}
From \cite{Fujiwara79}, \cite{Fujiwara}, there exists
$\delta >0$ independent of $\e$ 
such that for $|t|\leq \delta$, 
\begin{equation}\label{eq:solfond}
{\tt U}^\e(t)f(x)= e^{-in\frac{\pi}{4}\operatorname{sgn}t}
\frac{1}{\left|2\pi \e t \right|^{n/2}}\int_{\R^n}
k^\e(t,x,y) e^{iS(t,x,y)/\e}  f(y)dy \ ,
\end{equation}
where $S$ solves the eikonal equation
$$\d_t S +\frac{1}{2}|\nabla_x S|^2 +{\tt V}(t,x)=0,$$
and  $k^\e$ is bounded as well as all its $(x,y)$-derivatives,
uniformly for $\eps\in ]0,1]$ and $|t|\leq \delta$.

The group ${\tt U}^\e$ is unitary on
$L^2(\R^n)$, and there
exist $\delta >0$ and $C>0$ independent of $\e \in ]0,1]$ such that
for $|t|\leq \delta $, 
$$\|U^\e(t)\|_{L^1\to L^\infty} \leq \frac{C}{|\e t|^{n/2}}. $$
As noticed in \cite{Caz} (see also \cite{KT}), this yields Strichartz type
inequalities for ${\tt U}^\e$. 
\begin{defin}\label{def:adm}
 A pair $(q,r)$ is {\bf admissible} if $2\leq r
  <\frac{2n}{n-2}$ (resp. $2\leq r\leq \infty$ if $n=1$, $2\leq r<
  \infty$ if $n=2$)  
  and 
$$\frac{2}{q}=\delta(r)\equiv n\left( \frac{1}{2}-\frac{1}{r}\right).$$
\end{defin}
The following proposition is a consequence of \eqref{eq:solfond} and
\cite{KT}. 
\begin{prop}[Strichartz inequalities]\label{prop:strichartz} The group 
${\tt U}^\e(t)$ satisfies:\\
1. For any admissible pair $(q,r)$, any finite interval $I$,
  there exists $C_r(I)$ such that 
\begin{equation}\label{eq:strichlib}
   \eps^{\frac{1}{q}} \left\| {\tt U}^\eps(t)u\right\|_{L^q(I;L^r)}\leq C_r(I)
   \|u\|_{L^2}. 
  \end{equation}
2. For any admissible pairs $(q_1,r_1)$ and $(q_2,r_2)$, and any
   finite interval $I$, there exists $C_{r_1,r_2}(I)$ such that 
\begin{equation}\label{eq:strichnl}
      \eps^{\frac{1}{q_1}+\frac{1}{q_2}}\left\| \int_{I\cap\{s\leq
      t\}} {\tt U}^\eps(t-s)F(s)ds 
      \right\|_{L^{q_1}(I;L^{r_1})}\leq C_{r_1,r_2}(I) \left\|
      F\right\|_{L^{q'_2}(I;L^{r'_2})}. 
    \end{equation}
The above constants are independent of $\eps$. 
\end{prop}
\begin{rema}
In the case of $V$, the above constants do depend on the length of the
time interval $I$ as soon as $\delta_j=1$ for at least one integer
$j$. 
\end{rema}
For $(q,r)$ an admissible pair and $I$ a time interval, define 
\begin{equation*}
%\begin{split}
Y_{r}(I) := \Big\{ \psi \in C(I;\Sigma); \
B\psi  \in 
L^q(I;L^r)\cap L^\infty(I;L^2),\
 \forall B\in 
\{ Id,\nabla_x,|x|\}\Big\}\ .
%\end{split}
\end{equation*}
The main result of this section is the following. 
\begin{prop}\label{prop:localex}
Let ${\tt V}$ satisfying Assumption~\ref{hyp:gen},  $\varphi$ and
$\si$ satisfying Assumption~\ref{hyp:nl}. There exist $T>0$ and a
unique solution
$\psi \in Y_{2\sigma+2}(]-T,T[)$
to the initial value problem,
\begin{equation}\label{eq:sanseps}
\left\{
\begin{split}
i\d_t \psi +\frac{1}{2}\Delta \psi &= {\tt V}(t,x)\psi + |\psi|^{2\si}\psi\
,\\
\psi_{\mid t=0}&= \varphi\ .
\end{split}
\right.
\end{equation}
This solution actually belongs to $Y(]-T,T[)$, where
\begin{equation*}
%\begin{split}
Y(I) := \Big\{ \psi \in C(I;\Sigma); \
B\psi  \in  
L^q(I;L^r),\ 
 \forall B\in 
\{ Id,\nabla_x,|x|\}, \forall (q,r)\textrm{  admissible}\Big\}.
%\end{split}
\end{equation*}
If the potential ${\tt V}$ does not depend on time, we
have the following conservation laws:
\begin{itemize}
\item Mass: $\dis \|\psi(t)\|_{L^2}=\|\varphi\|_{L^2}$, $\forall
|t|<T$.
\item Energy: $$ E(t):=
\frac{1}{2}\|\nabla_x\psi(t)\|^2_{L^2}+\frac{1}{\si+1}\|\psi(t)\|_{L^{2\si
+2}}^{2\si
+2} +\int {\tt V}(x)|\psi(t,x)|^2dx \equiv E(0),\ \ \forall 
|t|<T.$$
\end{itemize}
\end{prop}
\begin{proof}
First, notice that Duhamel's principle for \eqref{eq:sanseps} writes
\begin{equation}\label{eq:duh1}
\psi(t,x)={\tt U}(t)\varphi - i\int_0^t{\tt U}(t-s)
\left(|\psi|^{2\si}\psi\right) (s)ds,
\end{equation}
where ${\tt U}(t):= {\tt U}^1(t)$. To estimate the nonlinear term, we
use Gagliardo-Nirenberg inequalities, which demand estimates on
$\nabla_x\psi$. We have,
\begin{equation*}
%\begin{split}
 \left[i \d_t +\frac{1}{2}\Delta -{\tt V}(t,x), \nabla_x
\right]=\nabla_x {\tt V}(t,x)\ ; \ \
 \left[i \d_t +\frac{1}{2}\Delta -{\tt V}(t,x), x
\right]=\nabla_x \ .
%\end{split}
\end{equation*}
Therefore, Duhamel's principles for $\nabla_x \psi$ and $x\psi$ are,
for $B\in \{ \nabla_x,x\}$,  
\begin{equation}\label{eq:duh2}
\begin{split}
B\psi(t,x)=&{\tt U}(t)B\varphi - i\int_0^t{\tt U}(t-s)
B\left(|\psi|^{2\si}\psi\right) (s)ds +i\int_0^t{\tt
U}(t-s)h_B(s)ds,\\
\textrm{with }&\ \  h_\nabla(t,x) = \nabla_x {\tt V}(t,x) \psi(t,x) \ ,\ \ 
h_x(t,x) = \nabla_x \psi(t,x). 
\end{split}
\end{equation}
Recall from Assumption~\ref{hyp:gen}, the potential ${\tt
V}$ is subquadratic, $\nabla_x {\tt V}(t,x) =O(\langle x \rangle )$,
locally in time. We formally have to solve a closed system of three equations
with three unknowns. This is achieved thanks to Strichartz
inequalities, provided by the case $\e =1$ in
Proposition~\ref{prop:strichartz}. The method is classical, and we
refer to \cite{Caz} for 
a complete proof. 
\end{proof}
\section{Inside the boundary layer}\label{sec:cl}

In this section, we prove that for any $\Lambda >0$, the solution
$u^\e$ to \eqref{eq:pb} is in
$C([-\Lambda\e,\Lambda\e];\Sigma)$ for $\e$ sufficiently small, and
satisfies the asymptotics \eqref{eq:asym1}.

Introduce the remainder $w^\e := u^\e -v^\e$. From
Proposition~\ref{prop:localex}, there exists $T^\e>0$ such that
$u^\e\in C([-T^\e,T^\e];\Sigma)$. Recall that $v^\e$ is given by
\eqref{eq:v}, where $\psi$ is the solution to  
\begin{equation}\label{eq:psi}\left\{
\begin{split}
i\d_t \psi +\frac{1}{2}\Delta \psi &= 
|\psi|^{2\sigma} \psi,\\
\psi_{\mid t=0} &= \varphi (x). 
\end{split}\right.
\end{equation}
It is well-known (see e.g. \cite{Caz}) that if $\varphi \in \Sigma$,
then $\psi \in C(\R,\Sigma)$, therefore $v^\e\in C(\R;\Sigma)$, and
$w^\e\in 
C([-T^\e,T^\e];\Sigma)$. 
This remainder solves
\begin{equation*}
\left\{
\begin{split}
i\e \d_t w^\e +\frac{1}{2}\e^2 \Delta w^\e &= V(x)u^\e
+ \e^{n\si} \left(|u^\e|^{2\si} u^\e-|v^\e|^{2\si} v^\e\right) \ ,\\
w^\e_{\mid t=0}&= 0\ .
\end{split}
\right.
\end{equation*}
We rewrite this problem as 
\begin{equation}\label{eq:pbw}
\left\{
\begin{split}
i\e \d_t w^\e +\frac{1}{2}\e^2 \Delta w^\e &= V(x)w^\e +V(x)v^\e
+ \e^{n\si} \left(|u^\e|^{2\si} u^\e-|v^\e|^{2\si} v^\e\right) \ ,\\
w^\e_{\mid t=0}&= 0\ .
\end{split}
\right.
\end{equation}
We shall actually prove a more precise result than that
stated in Theorem~\ref{theo}. 
\begin{prop}\label{prop:inside}
Suppose that Assumptions~\ref{hyp:pot} and \ref{hyp:nl} are
satisfied. Let $\Lambda >0$. Then for $0<\e\leq \e(\Lambda)$, $u^\e\in
C([-\Lambda \e, \Lambda\e];\Sigma)$ and 
\begin{equation*}
\limsup_{\e \to 0}\sup_{|t|\leq \Lambda \e}\left(
\left\|w^\e(t) \right\|_{L^2}+
\left\|A^\e_{\ell ,j}(t) w^\e \right\|_{L^2} \right)
   =0\ , \ \forall (\ell,j)\in \{1,2\}\times\{1,\ldots,n\}. 
\end{equation*}
\end{prop}

Recall that $U^\e(t)$ is the group associated to the linear part of
\eqref{eq:pb}, given by \eqref{eq:group}. It
satisfies Strichartz inequalities stated in
Proposition~\ref{prop:strichartz}. Duhamel's principle for
\eqref{eq:pbw} is
\begin{equation}\label{eq:duhw}
\begin{split}
w^\e(t) = &-i\e^{n\si-1} \int_0^t U^\e(t-s)\left(|u^\e|^{2\si}
u^\e-|v^\e|^{2\si} v^\e\right)(s)ds\\
& -i\e^{-1}\int_0^t
U^\e(t-s)V(x)v^\e(s)ds.
\end{split} 
\end{equation}
To apply the results of Proposition~\ref{prop:strichartz}, we
introduce special indexes in  the following algebraic lemma, whose
easy proof is left out. 
\begin{lem}\label{lem:alg}
Let $\si$ as in Assumption~\ref{hyp:nl}. 
There exist $\underline{q}$, $\underline{r}$, $\underline{s}$ and
$\underline{k}$ satisfying
\begin{equation}\label{eq:holder}
  \left\{
  \begin{split}
    \frac{1}{\underline{r}'}&=\frac{1}{\underline{r}}+
\frac{2\sigma}{\underline{s}},\\ 
    \frac{1}{\underline{q}'}&=\frac{1}{\underline{q}}+
\frac{2\sigma}{\underline{k}},
  \end{split}\right.
\end{equation}
and the additional conditions:
\begin{itemize}
\item The pair $(\underline{q},\underline{r})$ is admissible,
\item $0<\frac{1}{\underline{k}}<\delta(\underline{s})<1$. 
\end{itemize}
If $n=1$, we choose $(\underline{q},\underline{r})=(\infty ,2)$,
$\underline{s}= \infty$ and $\underline{k} = 2\sigma$.
\end{lem}
From Proposition~\ref{prop:strichartz} applied with the above indexes,
and H\"older inequality, \eqref{eq:duhw} yields, for $I^\e\ni 0$ a
time interval contained in $[-T^\e,T^\e]$, 
\begin{equation}\label{eq:estw1}
\begin{split}
\|w^\e\|_{L^{\underline q}(I^\e;L^{\underline r})}\lesssim
&\  \e^{n\si -1 -2/{\underline q}}
\left(\|u^\e\|_{L^{\underline k}(I^\e;L^{\underline s})}^{2\si} +
  \|v^\e\|_{L^{\underline k}(I^\e;L^{\underline s})}^{2\si} 
\right)\|w^\e \|_{L^{\underline q}(I^\e;L^{\underline r})}\\
& +
\e^{-1-1/{\underline q}}\|Vv^\e\|_{L^1(I^\e;L^2)} . 
\end{split}
\end{equation}
We now have two tasks:
\begin{itemize}
\item Estimate the source term $\|Vv^\e\|_{L^1(I^\e;L^2)}$.
\item Control the factor $ \|u^\e\|_{L^{\underline
k}(I^\e;L^{\underline s})}^{2\si} + 
  \|v^\e\|_{L^{\underline k}(I^\e;L^{\underline s})}^{2\si}$.  
\end{itemize}
Recall that $v^\e$ is given by \eqref{eq:v}, so 
\begin{equation*}
\left\| V(\cdot)v^\e(t,\cdot)\right\|^2_{L^2} =
\frac{1}{4}\sum_{j=1}^n  \om_j^4 \e^4 \left\| x_j^2 \psi (\e
t,x_j)\right\|^2_{L^2}+\sum_{j=1}^n b_j^2 \e^2 \left\|
x_j \psi (\e 
t,x_j)\right\|^2_{L^2} .  
\end{equation*}
If $\varphi \in \Sigma$, then $\psi \in C(\R,\Sigma)$, and the above
quantities are infinite in general. 

\subsection{Further regularity for $\psi$ when $\varphi \in {\EuScript
S}(\R^n)$} 
\label{sec:fur}
If we assume that $\varphi$ belongs to the Schwartz space
${\EuScript S}(\R^n)$, then we can prove additional regularity for
$\psi$. 
\begin{lem}\label{lem:fur}
Let $\varphi \in {\EuScript
S}(\R^n)$, and $\psi$ be the solution of the initial value problem
\eqref{eq:psi}.  Let $\si$ satisfying Assumption~\ref{hyp:nl}, and
$\Lambda >0$. Then, 
\begin{equation*}
\begin{split}
|x|^k \psi& \in
C([-\Lambda, \Lambda],L^2),\ \forall k\leq 3,\\
|x|^k \nabla_x \psi & \in
C([-\Lambda, \Lambda],L^2),\ \forall k\leq 2.
\end{split}
\end{equation*}
\end{lem}
\begin{proof}
As mentioned above, it is well-known that $\psi \in C([-\Lambda,
\Lambda],\Sigma)$. Using the simple remark,
$$\left[ i\d_t +\frac{1}{2}\Delta, x \right]= \nabla_x,$$
the function $x_j \psi$ solves, for $1\leq j \leq n$, 
\begin{equation}\label{eq:x_jpsi}
\left(i\d_t +\frac{1}{2}\Delta \right)x_j\psi = \d_j \psi
+|\psi|^{2\sigma} x_j \psi.
\end{equation}
For $1\leq k \leq n$, we have,
\begin{equation}\label{eq:x_jx_kpsi}
\left(i\d_t +\frac{1}{2}\Delta \right)x_jx_k\psi = \d_k( x_j \psi)+
x_k \d_j \psi
+|\psi|^{2\sigma} x_j x_k \psi.
\end{equation}
This shows that to know that $ x_jx_k\psi \in C([-\Lambda,
\Lambda],L^2)$, it is enough to prove that $x_\ell \nabla_x \psi \in
C([-\Lambda, 
\Lambda],L^2)$, for any $\ell$. Differentiating \eqref{eq:psi} with
respect to $x$ yields,
\begin{equation}\label{eq:nabla_xpsi}
\left(i\d_t +\frac{1}{2}\Delta \right)\nabla_x\psi = 
(\sigma +1)|\psi|^{2\sigma} \nabla_x \psi+ \sigma
|\psi|^{2\sigma-2}\psi^2 \overline{\nabla_x \psi} .
\end{equation}
Therefore,
\begin{equation}\label{eq:x_lnabla_xpsi}
\left(i\d_t +\frac{1}{2}\Delta \right)x_\ell\nabla_x\psi = \d_\ell
\nabla_x\psi +
(\sigma +1)|\psi|^{2\sigma} \nabla_x \psi+ \sigma
|\psi|^{2\sigma-2}\psi^2 \overline{\nabla_x \psi} .
\end{equation}
This shows that it is enough to know that $\Delta \psi \in C([-\Lambda,
\Lambda],L^2)$. This is well-known, from an idea due to Kato
(\cite{Kato87}, see also \cite{Caz}). The idea consists in differentiating
\eqref{eq:psi} with respect to time and proving that
$\d_t \psi \in C([-\Lambda,
\Lambda],L^2)$ when $\varphi \in H^2(\R^n)$. Then from
\eqref{eq:psi}, 
we deduce that $\Delta \psi \in C([-\Lambda,
\Lambda],L^2)$. Thus, $|x|^k \psi \in C([-\Lambda,
\Lambda],L^2)$ for $k \leq 2$ and $|x|^k \nabla_x \psi \in
C([-\Lambda, 
\Lambda],L^2)$ for $k \leq 1$.

Now, we can apply Kato's method to \eqref{eq:nabla_xpsi}, and prove that
if the
nonlinearity $F(z)=|z|^{2\sigma}z$ is twice differentiable (hence the
assumption $\sigma >1/2$ in Assumption~\ref{hyp:nl}), then $\d_t
\nabla_x \psi \in  C([-\Lambda, 
\Lambda],L^2)$. When using this information in \eqref{eq:x_jpsi},
Kato's method proves that 
$x_j \d_t \psi \in  C([-\Lambda,
\Lambda],L^2)$. Using the equation \eqref{eq:x_jpsi}, we deduce that
$x_j \Delta \psi \in  C([-\Lambda, \Lambda],L^2)$. This information is
enough to complete the proof of Lemma~\ref{lem:fur}. Multiplying
\eqref{eq:x_jx_kpsi} by $x_\ell$ yields,
\begin{equation*}\label{eq:x_jx_kx_lpsi}
\left(i\d_t +\frac{1}{2}\Delta \right)x_jx_kx_\ell\psi = 
\d_\ell(x_k \d_j \psi) + x_\ell \d_{k}( x_j \psi)+ x_\ell x_k \d_j \psi 
+|\psi|^{2\sigma} x_j x_k x_\ell \psi.
\end{equation*}
Reasoning as above, it is enough to know that $x \Delta \psi \in  C([-\Lambda,
\Lambda],L^2)$ and $x^\alpha \nabla_x \psi \in  C([-\Lambda,
\Lambda],L^2)$ for $|\alpha|\leq 2$. We saw how to prove the first
point. We know that the second holds for $|\alpha|\leq 1$, thus we just
have to multiply \eqref{eq:x_lnabla_xpsi} by $x_k$,
\begin{equation*}\label{eq:x_k x_lnabla_xpsi}
\begin{split}
\left(i\d_t +\frac{1}{2}\Delta \right)x_k x_l\nabla_x\psi = &\d_k
(x_l\nabla_x\psi) +
x_k\d_l
\nabla_x\psi \\
&+
(\sigma +1)|\psi|^{2\sigma} x_k\nabla_x \psi+ \sigma
|\psi|^{2\sigma-2}\psi^2 \overline{x_k\nabla_x \psi} .
\end{split}
\end{equation*}
Since $x \Delta \psi \in  C([-\Lambda,
\Lambda],L^2)$, we deduce that $|x|^2 \nabla_x \psi \in  C([-\Lambda,
\Lambda],L^2)$, which completes the proof.
\end{proof}
\begin{rema}
The assumption $\sigma >1/2$ could be removed if we considered a
smoother nonlinearity. Indeed, if we replaced $|\psi|^{2\sigma}\psi$
by $f(|\psi|^2)\psi$, with $f$ smooth and
$$f(|\psi|^2)\lesssim |\psi|^{2\sigma} \textrm{ when }|\psi|\rightarrow
0,$$
we could prove Lemma~\ref{lem:fur} without the assumption $\sigma
>1/2$, and even more regularity for $\psi$ (see for instance
\cite{HNT87}, \cite{HNT88}). This means, for \eqref{eq:pb}, that we
would replace $\eps^{n \sigma} |u^\eps|^{2\sigma}u^\eps$ by
$f(\eps^n |u^\eps|^{2})u^\eps$.
\end{rema}

We apply Lemma~\ref{lem:fur} to study \eqref{eq:pb} thanks to the following
result, which can be found for instance in \cite{Ginibre},
Proposition~3.5.
\begin{prop}\label{prop:gwp}
Let $\varphi$ and $\si$ satisfying Assumption~\ref{hyp:nl}. Let
$\delta >0$ and $\varphi_\delta \in {\EuScript S}(\R^n)$ such that 
$\Vert \varphi - \varphi_\delta \Vert_{\Sigma} \leq \delta$. If
$\psi_\delta$ denotes the solution to \eqref{eq:nls} with initial
datum $\varphi_\delta$, then 
\begin{equation*}
\left\Vert U_0(-t)\left(\psi(t) -\psi_\delta(t)\right)
\right\Vert_{L^\infty(\R;\Sigma
)}\Tend \delta 0 0\ ,
\end{equation*}
and in particular, for every $\Lambda >0$, 
\begin{equation*}
\Vert \psi -\psi_\delta\Vert_{L^\infty([-\Lambda,\Lambda];\Sigma
)}\Tend \delta 0 0\ .
\end{equation*}
\end{prop}

\subsection{The coupling term} \label{sec:in}
We want to estimate $ \|u^\e\|_{L^{\underline
k}(I^\e;L^{\underline s})}^{2\si} + 
  \|v^\e\|_{L^{\underline k}(I^\e;L^{\underline s})}^{2\si}$. \\ 
Gagliardo-Nirenberg inequalities and the
  property $\psi \in C(\R;\Sigma)$ yield
\begin{equation*}
\|v^\e(t)\|_{L^{\underline s}}= \e^{-\delta ({\underline s})}
\|\psi(\e t)\|_{L^{\underline s}}\lesssim \e^{-\delta ({\underline s})}
\|\psi(\e  t)\|_{L^2}^{1- \delta ({\underline s})} 
\|\nabla_x\psi(\e  t)\|_{L^2}^{\delta ({\underline s})} \leq C_\Lambda
  \e^{-\delta ({\underline s})}, 
\end{equation*}
for $|t|\leq \Lambda \e$, where $C_\Lambda$ does not depend on
$\e$. We expect a similar estimate to hold also for $u^\e$. From
\eqref{eq:sobolev}, it will be so if we know that $u^\e$ is bounded in
$L^2$, as well as $A^\e_{\ell,j}u^\e$ for any $\ell$ and $j$. The
first point is easy: so long as $u^\e$ is defined and sufficiently
smooth, its $L^2$-norm is constant (see
Proposition~\ref{prop:localex}). Showing the second is part of our
proof. Since $\psi \in C(\R;\Sigma)$, it is easy to check that for any
$\Lambda >0$, there exists $C(\Lambda)$ independent of $\e \in ]0,1]$,
such that 
\begin{equation*}
\|A^\e_{\ell,j}(t)v^\e(t,\cdot)\|_{L^2} \leq C(\Lambda)\ ,\ \ 
\forall (\ell,j)\in \{1,2\}\times\{1,\ldots,n\},\  \forall |t|\leq
\Lambda\e \ .
\end{equation*}
Since $w^\e=0$ at time $t=0$ and $w^\e \in C([-T^\e,T^\e];\Sigma)$ for
some $T^\e>0$, there exists $t^\e>0$ such that for $|t|<t^\e$, 
\begin{equation}\label{eq:tq1}
\|A^\e_{\ell,j}(t)w^\e(t,\cdot)\|_{L^2} \leq C(\Lambda)\ ,\ \ 
\forall (\ell,j)\in \{1,2\}\times\{1,\ldots,n\}\ .
\end{equation}
So long as \eqref{eq:tq1} holds, we can estimate
$\|u^\e(t)\|_{L^{\underline s}}$ like $\|v^\e(t)\|_{L^{\underline
s}}$, up to doubling the constants, but with the same power of $\e$. 

Let $\eta>0$ to be fixed later, and$I_\e \subset [-\eta \e,\eta\e]$
such that \eqref{eq:tq1} 
holds on $I_\e$. 
If $\varphi \in {\mathcal S}(\R^n)$, Lemma~\ref{lem:fur} and
\eqref{eq:estw1} yield
\begin{equation}\label{eq:estw2}
\|w^\e\|_{L^{\underline q}(I^\e;L^{\underline r})}\lesssim
\  \e^{n\si -1 -2/{\underline q}-2\si \delta ({\underline s})+
2\si/{\underline k}}\eta^{2\si/{\underline k}}
\|w^\e \|_{L^{\underline q}(I^\e;L^{\underline r})}\\
 +
\e^{1-1/{\underline q}}. 
\end{equation}
From Lemma~\ref{lem:alg}, 
$$n\si -1 -\frac{2}{{\underline q}}-2\si \delta ({\underline s})+
\frac{2\si}{{\underline k}}=0\ ,$$
and for $\eta >0$ sufficiently small, the first term of the right hand
side of \eqref{eq:estw2} is absorbed by the left hand side,
\begin{equation}\label{eq:estw3}
\|w^\e\|_{L^{\underline q}(I^\e;L^{\underline r})}\lesssim
\e^{1-1/{\underline q}}. 
\end{equation}
Apply Strichartz inequality \eqref{eq:strichnl} again, with now
$r_1=2$ and $r_2=\underline r$,
\begin{equation}\label{eq:estw4}
\|w^\e\|_{L^\infty(I^\e;L^2)}\lesssim
\e^{n\si -1 -1/{\underline q}-2\si \delta ({\underline s})+
2\si/{\underline k}}\eta^{2\si/{\underline k}}
\|w^\e \|_{L^{\underline q}(I^\e;L^{\underline r})}+ \e\lesssim \e\ , 
\end{equation}
from \eqref{eq:estw3}. 

Assuming for a moment that we know that \eqref{eq:tq1} holds for $|t|\leq
\Lambda\e$, the above computation, repeated a finite number of times,
 yields an estimate of the form 
\begin{equation}\label{eq:finw}
\|w^\e\|_{L^\infty([-\Lambda \e,\Lambda \e];L^2)}\leq C\e\, 
e^{C\Lambda}\ . 
\end{equation}

To prove that indeed \eqref{eq:tq1} holds for $|t|\leq
\Lambda\e$, we follow the same lines as above, replacing $w^\e$ by
$A^\e_{\ell ,j}w^\e$. Since $A^\e_{\ell ,j}$ commute with the linear
part of \eqref{eq:pb} (see the first point of Lemma~\ref{lem:ab}),
the analog of \eqref{eq:duhw} for $A^\e_{\ell ,j}w^\e$ is
\begin{equation}\label{eq:duhAw}
\begin{split}
A^\e_{\ell ,j}(t)w^\e = &-i\e^{n\si-1} \int_0^t U^\e(t-s)A^\e_{\ell
,j}(s)\left(|u^\e|^{2\si} 
u^\e-|v^\e|^{2\si} v^\e\right)(s)ds\\
& -i\e^{-1}\int_0^t
U^\e(t-s)A^\e_{\ell ,j}(s)\left(V(x)v^\e(s)\right)ds.
\end{split} 
\end{equation}
From Lemma~\ref{lem:fur}, the source term (the last term in the above
expression) is estimated as before. From \eqref{eq:jauge} and
\eqref{eq:sobolev}, we can estimate the first term of the right hand
side of \eqref{eq:duhAw} as above. This yields finally, so long as
\eqref{eq:tq1} holds and for $|t|\leq \Lambda \e$, 
\begin{equation}\label{eq:finAw}
\|A^\e_{\ell ,j} w^\e\|_{L^\infty(I^\e;L^2)}\leq
C(\ell ,j)\e\, 
e^{C(\ell, j)\Lambda}\ . 
\end{equation}

\subsection{Conclusion} Let $\delta >0$, and $\varphi_\delta \in
{\EuScript S}(\R^n)$ such that  
$\Vert \varphi - \varphi_\delta \Vert_{\Sigma} \leq \delta$. Define
$\psi_\delta$ as the solution to \eqref{eq:nls} with initial datum
$\varphi_\delta$, and $v^\e_\delta$ by 
$$v^\e_\delta (t,x)=\frac{1}{\e^{n/2}}\psi_\delta \left(\frac{t}{\e},
\frac{x}{\e}  \right)\ .$$
The remainder $w^\e_\delta := u^\e -v^\e_\delta$ solves
\begin{equation}\label{eq:pbwdelta}
\left\{
\begin{split}
i\e \d_t w^\e_\delta +\frac{1}{2}\e^2 \Delta w^\e_\delta &= V(x)w^\e_\delta
+V(x)v^\e_\delta 
+ \e^{n\si} \left(|u^\e|^{2\si} u^\e-|v^\e_\delta|^{2\si}
v^\e_\delta\right) \ ,\\ 
w^\e_{\delta}(0,x) &= \frac{1}{\e^{n/2}}(\varphi -
\varphi_\delta) \left(
\frac{x}{\e}  \right)\ ,
\end{split}
\right.
\end{equation} 
which is the analog of\eqref{eq:pbw}, with an initial datum which is
nonzero, but arbitrarily small in $\Sigma$ (as $\delta$ goes to
zero). 

Our method proves both the existence of $u^\e$ in $\Sigma$ up to time
$\Lambda \e$ for $\e$
sufficiently small, and the asymptotics \eqref{eq:asym1}. This
approach is classical in geometrical optics (see
e.g. \cite{RauchUtah}). 
From Proposition~\ref{prop:localex}, it is well defined in $\Sigma$ on
the time interval
$[-T^\e,T^\e]$ for some $T^\e >0$.
Since $v^\e_\delta \in C(\R;\Sigma)$, we want to prove that
$w^\e_\delta$ exists 
in $\Sigma$ up to time
$\Lambda \e$ for $\e$
sufficiently small, and is asymptotically small. By construction, we have
\begin{equation}\label{eq:large0}
\left\| w^\e_\delta(0) \right\|_{L^2} +\sum_{(\ell,j)\in
\{1,2\}\times\{1,\ldots,n\}} 
\left\|A^\eps_{\ell,j}(0)
w^\e_\delta \right\|_{L^2}\leq \delta\ .
\end{equation}
From Proposition~\ref{prop:localex}, either $ w^\e_\delta$ (hence
$u^\e$) exists in $\Sigma$ on the time interval $[-\Lambda \e,\Lambda
\e]$, or the 
maximal solution belongs to $C([0,T^\e[;\Sigma)$ with $0<T^\e<\Lambda \e$ and
\begin{equation*}
\liminf_{t \to T^\e} \| w^\e_\delta(t)\|_\Sigma =\infty\ .
\end{equation*}
In the latter case, for any $\Gamma >0$, there is a first time,
$T_\Gamma^\e$ such that 
\begin{equation}\label{eq:pv0}
\left\|\widetilde w^\e(T_\Gamma^\e
) \right\|_{L^2} +\sum_{(\ell,j)\in \{1,2\}\times\{1,\ldots,n\}}
\left\|A^\eps_{\ell,j}(T_\Gamma^\e)
\widetilde w^\e\right\|_{L^2}= \Gamma\delta\ .
\end{equation}
We prove that
there is $\Gamma >0$ independent of $\Lambda$ and $\e$, and a constant
$C=C(\Lambda)$ independent of $\e$ such that 
for $\e \leq 1$ and $t^\e\leq T_\Gamma^\e$, 
\begin{equation}\label{eqw}
\sup_{|t|\leq \Lambda \e}\left(\left\| w^\e_\delta(
t) \right\|_{L^2} +\sum_{(\ell,j)\in \{1,2\}\times\{1,\ldots,n\}}
\left\|A^\eps_{\ell,j}(t)
w^\e_\delta \right\|_{L^2}\right)\leq \frac{\Gamma}{2}\delta + C\e \ .
\end{equation}
Choosing $\e$ sufficiently small so that $C\e<\Gamma/2$ 
contradicts \eqref{eq:pv0}. This proves that we can take $t^\e
=\Lambda \e$ in
\eqref{eq:tq1}.

Resuming the computations of Section~\ref{sec:in} yields the same
estimates as \eqref{eq:estw2}, plus a term estimated by $\delta
\e^{-1/\underline q}$, due to the initial datum. This means that in
\eqref{eq:estw4}, \eqref{eq:finw} and \eqref{eq:finAw}, we have to
replace $\e$ by $\e 
+\delta$ in the right hand sides; this yields \eqref{eqw}. 
We infer, 
\begin{equation*}
\limsup_{\e \to 0}\sup_{|t|\leq \Lambda \e}\left(
\left\|w^\e_\delta(t) \right\|_{L^2}+
\sum_{(\ell,j)\in \{1,2\}\times\{1,\ldots,n\}}
\left\|A^\e_{\ell ,j}(t) w^\e_\delta \right\|_{L^2} \right)
   \leq \frac{\Gamma}{2}\delta \ ,  
\end{equation*}
where $C$ does not depend on $\delta$. Choosing $\delta$ arbitrarily
small, the above estimate and Proposition~\ref{prop:gwp} yield
Proposition~\ref{prop:inside}.

Finally, Proposition~\ref{prop:inside} implies the asymptotics
\eqref{eq:asym1}. Rewrite the definition of $A^\e_{\ell ,j}$,
\begin{equation*}
\left(
	\begin{array}{c}
	A^\e_{1,j}\\
	A^\e_{2,j}
	\end{array}\right)= 
\left(
	\begin{array}{cc}
	h_j&g_j/\e\\
	-\e\delta_j\om_j^2g_j& h_j
	\end{array}\right)\left(
	\begin{array}{c}
	x_j/\e\\
	i\e \d_j
	\end{array}\right) -b_j\left(
	\begin{array}{c}
	t^2/(2\e)\\
	t
	\end{array}\right).
\end{equation*}
The determinant of the above matrix is 
$$h_j^2 +\delta_j \om_j^2 g_j^2 \equiv 1,$$
and we have
\begin{equation}\label{eq:backop}
\begin{split}
\frac{x_j}{\e}& =
h_j(t)A_{1,j}^\e(t)-\frac{g_j(t)}{\e}A_{2,j}^\e(t)+b_j\left(
\frac{t^2}{2\e} h_j(t) -\frac{t}{\e}g_j(t)\right)\ ,\\
i\e \d_j &= \e \delta_j \om_j^2 g_j(t)A_{1,j}^\e(t) +
h_j(t)A_{2,j}^\e(t) +b_j\left(\delta_j \om_j^2 \frac{t^2}{2}g_j(t) +th_j(t)
\right)\ .
\end{split}
\end{equation}
Since $g_j(t)=O(t)$ as $t$ goes to zero, 
 it is clear that Proposition~\ref{prop:inside} implies the asymptotics
\eqref{eq:asym1}. 

\section{Beyond the boundary layer}\label{sec:beyond}

In this section, we complete the proof of Theorem~\ref{theo}. The end
of the proof is divided into two parts; we first study the transition
between the two r\'egimes \eqref{eq:asym1} and \eqref{eq:asym2}, then 
prove the existence of $u^\e$ along with the asymptotics
\eqref{eq:asym2}. Since the proofs are similar for positive or
negative times, we restrict to the case of positive times. 

\subsection{Matching the two r\'egimes}
\label{sec:matching}
In Proposition~\ref{prop:inside}, $\Lambda$ was a fixed parameter; in
any boundary layer of size $\Lambda\e$ around the origin, the
asymptotic behaviour of $u^\e$ is given by $v^\e$. For $t\gg \e$, the
behaviour of $u^\e$ is asymptotically the same as that of $u_+^\e$. We
now prove that the transition between these two r\'egimes occurs in
a boundary layer of size $\Lambda\e$, when $\Lambda$ goes to
infinity. 
\begin{prop}\label{prop:matching}
The function $u^\eps_+$ becomes an approximate solution of $u^\eps$
when $t$ reaches $\Lambda \eps$, for large $\Lambda$.
\begin{equation*}
\begin{split}
\limsup_{\eps \rightarrow 0} &\left( \left\|u^\eps(\Lambda
\eps) -u^\eps_+
(\Lambda \eps)\right\|_{L^2} +\left\|A^\eps_{\ell,j}(\Lambda \eps)
\left( u^\eps -u^\eps_+
\right)\right\|_{L^2} \right)\Tend \Lambda {+\infty} 0\ ,\\
& \forall
(\ell,j)\in \{1,2\}\times\{1,\ldots,n\}.  
\end{split}
\end{equation*}
\end{prop}
\begin{proof}
From Proposition~\ref{prop:inside}, we only have to prove the above
limit when $u^\eps$ is replaced by $v^\eps$. We proceed to another
reduction of the problem, by noticing that for $|t|\leq \Lambda \e$,
the role of the potential $V$ is negligible not only for $u^\e$, but
also for $u_+^\e$. Define $v_+^\e$ by 
\begin{equation}\label{eq:upapp}\left\{
\begin{split}
i\eps \d_t v_+^\eps +\frac{1}{2}\eps^2\Delta
v_+^\eps &=0,  \\
v^\eps_{+ \mid t=0} &= \frac{1}{\eps^{n/2}}\psi_+
\left( 
\frac{x}{\eps}\right).
\end{split}\right.
\end{equation}
By scaling, we have 
$$v_+^\eps (t,x) =\frac{1}{\eps^{n/2}}\psi_+^0
\left(\frac{t}{\eps}, \frac{x}{\eps}\right),$$
where $\psi_+^0(t,x)=\exp(it\Delta/2)\psi_+(x)$. 
\begin{lem}\label{lem:cl}
Let $\Lambda \geq 1$. 
The potential $V$ is negligible for $0\leq t\leq \Lambda \eps$ in
\eqref{eq:upm},
\begin{equation*}
\begin{split}
\limsup_{\eps \rightarrow 0}\sup_{0\leq t\leq \Lambda \eps}
 &\left( \left\|u^\eps_+(t) -v^\eps_+
(t)\right\|_{L^2} +\left\|A^\eps_{\ell,j}(t)
\left( u^\eps_+ -v^\eps_+
\right)\right\|_{L^2} \right)=0\ ,\\
& \forall
(\ell,j)\in \{1,2\}\times\{1,\ldots,n\}.  
\end{split}
\end{equation*}
\end{lem}
\begin{proof}[Proof of Lemma~\ref{lem:cl}]
Denote $w_+^\eps = u_+^\eps - v_{+}^\eps$. We have,
\begin{equation*}\left\{
\begin{split}
i\eps \d_t  w_+^\eps+\frac{1}{2}\eps^2\Delta
w_+^\eps &=V(x)w_+^\eps + V(x)v_{+}^\eps,  \\
w^\eps_{+ \mid t=0} &= 0.
\end{split}\right.
\end{equation*}
From the classical energy estimates (which are also a consequence of 
Strichartz inequalities), 
\begin{equation*}
\begin{split}
\sup_{0\leq t\leq \Lambda \eps}\|w_+^\eps(t)\|_{L^2} & \lesssim 
\eps^{-1}\int_0^{\Lambda
\eps}\|V(.)v_{+\rm{app}}^\eps(t,.)\|_{L^2}dt\\
& \lesssim \int_0^{\Lambda }\|V(\eps \cdot)\psi_+^0(t,\cdot
)\|_{L^2}dt\ . 
\end{split}
\end{equation*}
By density (for $\psi_+$), we can assume that $\psi^0_+$ has the same
smoothness as in Lemma~\ref{lem:fur} (the proof is even easier since
we now consider \emph{linear} problems). In that case we have
$$\sup_{0\leq t\leq \Lambda \eps}\|w_+^\eps(t)\|_{L^2} =O(\e)\ .$$
The proof that $A^\eps_{\ell,j}(t)w_+^\eps$ satisfies the same
property is straightforward. Finally, without the smoothness
assumption of Lemma~\ref{lem:fur}, $O(\e)$ is replaced by $o(1)$, and
the proof of Lemma~\ref{lem:cl} is complete. 
\end{proof}
%To complete the proof of Proposition~\ref{prop:matching}, we use an
%asymptotic property of the free Schr\"odinger group
%$\exp(it\Delta/2)$, for the proof of which we refer to \cite{Rauch91}
%for instance. 
%\begin{lem}\label{lem:LS}
%Let $\varphi \in \Sigma$, and denote $U_0(t):= \exp(it\Delta/2)$. The
%following asymptotics holds, 
%$$\left\|\varphi  -U_0(-t)\left(\frac{e^{-i\frac{x^2}{2t}}}{t^{n/2}} \widehat
%\varphi\left( \frac{x}{t}\right)\right)\right\|_{\Sigma}\Tend t {+
%\infty} 0\ ,$$ 
%where the Fourier transform is defined by 
%\begin{equation*}
%{\EuScript F} f (\xi) =\widehat f(\xi)=
%  \frac{1}{(2i\pi)^{n/2}}\int_{\R^n} e^{-ix.\xi}f(x)dx.
%\end{equation*}
%\end{lem}
Recall that we have
\begin{equation*}
\begin{split}
v_{+}^\eps(\Lambda \eps,x)&=\frac{1}{\eps^{n/2}}U_0(\Lambda
)\psi_+\left( \frac{x}{\eps}\right)\ ,\
v^\e(\Lambda \eps,x)=\frac{1}{\eps^{n/2}}\psi \left(\Lambda,
\frac{x}{\e}\right)\ ,\\
&\lim_{t\to +\infty}\Big\| U_0(-t)\psi(t)-\psi_+ \Big\|_\Sigma
=0 \ ,
\end{split}
\end{equation*}
where the last line is nothing but \eqref{eq:scatt}. This implies in
particular, since $U_0$ is unitary on $L^2$,
\begin{equation*}
\limsup_{\e \to 0}\left\| v^\e(\Lambda \e)-v^\e_+(\Lambda
\e)\right\|_{L^2} \Tend \Lambda {+\infty} 0\ ,
\end{equation*}
which is the first asymptotics in Proposition~\ref{prop:matching}. 

To conclude the proof, the idea is that the operator appearing in
\eqref{eq:scatt} are close to the operators $A^\e_{\ell,j}(t)$ for
$|t|\leq \Lambda \e$. Using the identity 
$$U_0(t)xU_0(-t)= x+it \nabla_x,$$
and the fact that the group $U_0$ is unitary on $L^2$, we can rewrite
\eqref{eq:scatt} as
\begin{equation*}
\begin{split}
\|\psi(t) -U_0(t)\psi_+\|_{L^2}&+ \|\nabla_x\psi(t) -U_0(t)\nabla_x
\psi_+\|_{L^2}\\
&+
\left\|\left(x+it\nabla_x \right)\left(\psi(t) -U_0(t)\psi_+\right)
\right\|_{L^2}\Tend t {+\infty} 0.
\end{split} 
\end{equation*}
From the definition of the function $h_j$'s and $g_j$'s, we have, as
$t\to 0$, 
$$h_j(t)=1 +O(t)\  \ \ ;\ \ g_j(t)=t +O(t^2) \ .$$
Therefore, we have, in the case of $v^e$,
\begin{equation*}
\begin{split}
\Big\|&\left(A_{1,j}^\e(\Lambda \e)-\frac{x_j}{\e}-i(\Lambda \e)\d_j
\right)
v^e(\Lambda\e,\cdot)\Big\|_{L^2}=\\
&= \left\|\left(\frac{x_j}{\e}(h_j(\Lambda \e) -1) +
i(g_j(\Lambda\e)-\Lambda\e)\d_j -\frac{b_j}{2}\Lambda^2\e
\right)v^e(\Lambda\e,\cdot)\right\|_{L^2} \\
&= \left\|\left(x_j(h_j(\Lambda \e) -1) +
i\frac{g_j(\Lambda\e)-\Lambda\e}{\e}\d_j -\frac{b_j}{2}\Lambda^2\e
\right)\psi(\Lambda,\cdot)\right\|_{L^2}\\
&= O(\e),
\end{split}
\end{equation*}
for any fixed $\Lambda \geq 1$, since $\psi \in C(\R;\Sigma)$. Similar
computations hold with $A_{2,j}^\e$, and when $v^e$ is replaced by
$v^e_+$. The proof of Proposition~\ref{prop:matching} is complete. 
\end{proof}

\subsection{The linear r\'egime} We now complete the proof of
Theorem~\ref{theo}. Fix $T>0$. From \eqref{eq:backop}, it is
enough to prove that $u^\e(t)$, as well as $A_{\ell,j}^\e(t)u^\e$
for any $\ell,j$, remains bounded in $L^2$, up to time $T$, provided
that $\e$ is sufficiently small. The relation \eqref{eq:backop} shows
in addition that we can prove the asymptotics \eqref{eq:asym2}
when the operators $\e\nabla$ and $x$ are replaced by the
$A_{\ell,j}^\e(t)$'s. 

Our method is the same as in Section~\ref{sec:cl}. 
Introduce the remainder 
$$\widetilde w^\e_+= u^\e -u_+^\e.$$ 
From Proposition~\ref{prop:inside}, it is well defined in $\Sigma$ up
to time $\Lambda \e$ for any $\Lambda >0$, provided that $\e$ is
sufficiently small. It solves
\begin{equation*}
i\e \d_t \widetilde w^\e +\frac{1}{2}\e^2 \Delta \widetilde w^\e =
V(x)\widetilde w^\e
+ \e^{n\si}|u^\e|^{2\si} u^\e \ . 
\end{equation*}
Since $v^\e_+ \in C(\R;\Sigma)$ (see in particular \eqref{eq:commut}
and \eqref{eq:backop}), we want to prove that $\widetilde w^\e$ exists
in $\Sigma$ up to time
$T$ for $\e$
sufficiently small, and is asymptotically small in the sense of
\eqref{eq:asym2}. From Proposition~\ref{prop:matching}, 
\begin{equation*}
%\begin{split}
\limsup_{\eps \rightarrow 0} \left( \left\|\widetilde w^\e(\Lambda
\eps) \right\|_{L^2} +\left\|A^\eps_{\ell,j}(\Lambda \eps)
\widetilde w^\e\right\|_{L^2} \right)\Tend \Lambda {+\infty} 0\ ,\
 \forall
(\ell,j)\in \{1,2\}\times\{1,\ldots,n\}.  
%\end{split}
\end{equation*}
Let $\delta >0$. From Proposition~\ref{prop:matching}, there exist
$\e_0>0$ and $\Lambda_0$ such that for  $0<\e\leq \e_0$ and
$\Lambda\geq \Lambda_0$, 
\begin{equation}\label{eq:large}
\left\|\widetilde w^\e(\Lambda
\eps) \right\|_{L^2} +\sum_{(\ell,j)\in \{1,2\}\times\{1,\ldots,n\}}
\left\|A^\eps_{\ell,j}(\Lambda \eps)
\widetilde w^\e\right\|_{L^2}\leq \delta\ .
\end{equation}
From Proposition~\ref{prop:matching} again, there exists $t^\e >
\Lambda \eps$ such that 
\begin{equation}\label{eq:tantque2}
\sup_{\Lambda \e \leq t\leq t^e}\left(\left\|\widetilde w^\e(
t) \right\|_{L^2} +\sum_{(\ell,j)\in \{1,2\}\times\{1,\ldots,n\}}
\left\|A^\eps_{\ell,j}(t)
\widetilde w^\e\right\|_{L^2}\right)\leq 2\delta\ .
\end{equation}
Let $0<\e\leq \e_0$ and
$\Lambda\geq \Lambda_0$. 
From Proposition~\ref{prop:localex}, either $\widetilde w^\e$ (hence
$u^\e$) exists in $\Sigma$ on the time interval $[0,T]$, or the
maximal solution belongs to $C([0,T^\e[;\Sigma)$ with $0<T^\e<T$ and
\begin{equation*}
\liminf_{t \to T^\e} \|\widetilde w^\e(t)\|_\Sigma =\infty\ .
\end{equation*}
From \eqref{eq:backop}, in the latter case, there is a first time,
$T_0^\e$ such that 
\begin{equation}\label{eq:pv}
\left\|\widetilde w^\e(T_0^\e
) \right\|_{L^2} +\sum_{(\ell,j)\in \{1,2\}\times\{1,\ldots,n\}}
\left\|A^\eps_{\ell,j}(T_0^\e)
\widetilde w^\e\right\|_{L^2}= 4\delta\ .
\end{equation}
We prove that, up to choosing $\Lambda$ even larger, 
there is a constant $C=C(T)$ independent of $\e$ and $\Lambda$ such that
for $\e \leq 1$ and $t^\e\leq T_0^\e$, 
\begin{equation}\label{eq:lk}
\sup_{\Lambda \e \leq t\leq t^\e}\left(\left\|\widetilde w^\e(
t) \right\|_{L^2} +\sum_{(\ell,j)\in \{1,2\}\times\{1,\ldots,n\}}
\left\|A^\eps_{\ell,j}(t)
\widetilde w^\e\right\|_{L^2}\right)\leq 3\delta + C\e^{2\si
\delta(\underline s)/n} \ .
\end{equation}
Choosing $\e$ sufficiently small so that $C\e^{2\si
\delta(\underline s)/n} <\delta$ 
contradicts \eqref{eq:pv}. This proves that we can take $t^\e =T$ in
\eqref{eq:tantque2},
hence the first point of
Theorem~\ref{theo}, along with the asymptotics \eqref{eq:asym2}, since
$\delta>0$ is arbitrary (recall that for any fixed $\delta>0$, we have
to choose $\e$ small and $\Lambda$ large, so that \eqref{eq:large}
holds).

Recall that $u^\e_+$ solves the linear equation \eqref{eq:upm}; its
$L^2$-norm is independent of time, and from \eqref{eq:commut}, the
same holds for $A^\e_{\ell,j}u^\e_+$, for any $\ell$ and $j$. 
So long as \eqref{eq:tantque2} holds, we thus have an $L^2$ bound for
$u^\e$ and $A^\e_{\ell,j}u^\e$,
\begin{equation}\label{eq:borneu}
\sup_{\Lambda \e \leq t\leq t^e}\left(\left\|u^\e(
t) \right\|_{L^2} +\sum_{(\ell,j)\in \{1,2\}\times\{1,\ldots,n\}}
\left\|A^\eps_{\ell,j}(t)
u^\e\right\|_{L^2}\right)\leq C_*\ .
\end{equation} 
Denote $J^\e:= [\Lambda\e,t^e]$.
From Strichartz inequalities and Lemma~\ref{lem:alg}, 
\begin{equation*}
\|\widetilde w^\e\|_{L^\infty(J^\e;L^2)}
\leq \|\widetilde w^\e(\Lambda \e)\|_{L^2} +
C\e^{n\si -1-1/{\underline q}}
\|u^\e\|_{L^{\underline k}(J^\e;L^{\underline s})}^{2\si} 
\|u^\e \|_{L^{\underline q}(J^\e;L^{\underline r})}\ . 
\end{equation*}
From \eqref{eq:borneu}, \eqref{eq:sobolev} and Lemma~\ref{lem:P}, we
infer that if $t^\e\leq \pi/(2\underline \om)$, 
\begin{equation*}
\|\widetilde w^\e\|_{L^\infty(J^\e;L^2)}
\leq \|\widetilde w^\e(\Lambda \e)\|_{L^2} +
\rho(\Lambda)\e^{1/\underline q}\|u^\e \|_{L^{\underline
q}(J^\e;L^{\underline r})}\ ,  
\end{equation*}
where $\rho(\Lambda)$ is a function independent
of $\e$ that goes to zero as $\Lambda$ goes to infinity. Using
\eqref{eq:borneu} and \eqref{eq:sobolev} again, we have
\begin{equation*}
\e^{1/\underline q}\|u^\e \|_{L^{\underline
q}(J^\e;L^{\underline r})}\leq C |J^\e|^{1/\underline q}\leq C
T^{1/\underline q} \ .
\end{equation*}
Therefore, 
\begin{equation*}
\|\widetilde w^\e\|_{L^\infty(J^\e;L^2)}
\leq \|\widetilde w^\e(\Lambda \e)\|_{L^2} +
CT^{1/\underline q} \rho(\Lambda)\ .
\end{equation*}
Taking
$\Lambda$ even larger if necessary, \eqref{eq:large} implies that if
$t^\e\leq \pi/(2\underline \om)$, then 
\begin{equation*}
\|\widetilde w^\e\|_{L^\infty(J^\e;L^2)}
\leq \|\widetilde w^\e(\Lambda \e)\|_{L^2}+\delta \ .
\end{equation*}
For $t^\e\geq
\pi/(2\underline \om)$, the second part of
Lemma~\ref{lem:P} implies  
\begin{equation}\label{qwer}
\begin{split}
\|\widetilde w^\e\|_{L^\infty(J^\e;L^2)}
 \leq \|\widetilde w^\e(\Lambda \e)\|_{L^2} +
& C\e^{n\si -1-1/{\underline q}}
\|u^\e\|_{L^{\underline k}(J^\e;L^{\underline s})}^{2\si} 
\|u^\e \|_{L^{\underline q}(J^\e;L^{\underline r})}\\
 \leq \|\widetilde w^\e(\Lambda \e)\|_{L^2} + & CT^{1/\underline q}
\e^{n\si -1-2/{\underline q}} 
\|u^\e\|_{L^{\underline k}(J^\e;L^{\underline s})}^{2\si} \\
 \leq \|\widetilde w^\e(\Lambda \e)\|_{L^2} + & CT^{1/\underline q}
\e^{n\si -1-2/{\underline q}} 
\|u^\e\|_{L^{\underline k}([\Lambda \e,\pi/(2\underline \om)]
;L^{\underline s})}^{2\si}\\
 +& CT^{1/\underline q}\e^{n\si -1-2/{\underline q}} 
\|u^\e\|_{L^{\underline k}([\pi/(2\underline \om),t^\e]
;L^{\underline s})}^{2\si}  \\
 \leq \|\widetilde w^\e(\Lambda \e)\|_{L^2} + &\delta +  C(T)\e^{2\si
\delta(\underline s)/n}\ .
\end{split}
\end{equation}
Computations for $A_{\ell,j}^\e (t)\widetilde w^\e$ are similar. Since
$A_{\ell,j}^\e$ acts like a derivative on the nonlinear term
(Lemma~\ref{lem:ab}), we have
\begin{equation*}
\|A_{\ell,j}^\e\widetilde w^\e\|_{L^\infty(J^\e;L^2)}
\leq \|A_{\ell,j}^\e(\Lambda \e)\widetilde w^\e\|_{L^2} +
C\e^{n\si -1-1/{\underline q}}
\|u^\e\|_{L^{\underline k}(J^\e;L^{\underline s})}^{2\si} 
\|A_{\ell,j}^\e u^\e \|_{L^{\underline q}(J^\e;L^{\underline r})}\ .  
\end{equation*}
Estimate \eqref{eq:borneu}, along with Proposition~\ref{prop:localex},
implies that there exists $C(T)$ such that for $t^\e \leq T$, 
\begin{equation*}
\e^{1/\underline q}\|A_{\ell,j}^\e u^\e \|_{L^{\underline
q}(J^\e;L^{\underline r})}\leq C(T)\ .
\end{equation*}
We thus have the same estimate as above, for $\Lambda$ sufficiently
large, 
\begin{equation}\label{asd}
\|A_{\ell,j}^\e\widetilde w^\e\|_{L^\infty(J^\e;L^2)}
\leq \|A_{\ell,j}^\e(\Lambda \e)\widetilde w^\e\|_{L^2} +
\frac{\delta}{2n} + C(T)\e^{2\si
\delta(\underline s)/n}\ .
\end{equation}
Summing \eqref{qwer} and \eqref{asd} yields \eqref{eq:lk}, which
completes the proof of Theorem~\ref{theo}.

\section{Partial results for general subquadratic potentials}
\label{sec:general}

Intuitively, there is no reason why Theorem~\ref{theo} should not be
true for more general potentials than \eqref{def:V1}, in particular for
potentials satisfying Assumption~\ref{hyp:gen}. We prove in particular
that \eqref{eq:asym1} still holds for this class of
potentials. However, we cannot prove  \eqref{eq:asym2}. From the
technical point of view, this is due to the lack of operators such as
$A_{\ell,j}^\e$. For the linear r\'egime, these operators have three
major advantages:
\begin{itemize}
\item They commute with the linear part of the equation, including the
potential, see \eqref{eq:commut}.
\item They yield modified Gagliardo-Nirenberg inequalities,
\eqref{eq:sobolev}. 
\item They act on the nonlinear term like derivatives, \eqref{eq:jauge}.
\end{itemize}
As we mentioned in the proof of Lemma~\ref{lem:ab}, the last two
points follow from the formula \eqref{eq:factor}. We first prove that
there exists an operator satisfying a similar formula \emph{and}
commuting with the linear part of the equation, \eqref{eq:commut}, if
and only if the potential is of the form we consider,
\eqref{def:V1}. We then prove \eqref{eq:asym1} for general potentials
satisfying Assumption~\ref{hyp:gen}.

\subsection{Lemma~\ref{lem:ab} holds only for potentials of the form
\eqref{def:V1}}\label{sec:eik}
 Let ${\tt V}$ satisfying
Assumption~\ref{hyp:gen}, \emph{independent of time} (${\tt V}={\tt
V}(x)$), and define an operator $A^\e(t)$ by 
\begin{equation}\label{eq:opgen}
A^\e(t) = i f(t)e^{i\phi(t,x)/\e}\nabla_x \left(e^{-i\phi(t,x)/\e} \cdot
\right)=\frac{f(t)}{\e}\nabla_x \phi(t,x) + i f(t)\nabla_x \ , 
\end{equation}
where $f$ and $\phi$ are real-valued functions, to be determined. This
is the generalization 
of  \eqref{eq:factor}. Such an operator formally satisfies
\eqref{eq:jauge} and an analog to \eqref{eq:sobolev}. 
Notice that in \eqref{eq:factor}, the phases $\phi_\ell$ ($\ell=1$ or
$2$) solve the eikonal equation
\begin{equation}\label{eq:eikon}
\d_t \phi +\frac{1}{2}|\nabla_x \phi|^2 +{\tt V}(x)=0\ .
\end{equation} 
\begin{prop}
Let $\phi\in C^4(]0,T]\times \R^n;\R)$ and $f\in C^1(]0,T])$ for some
$T>0$. Assume that $f$ does not cancel on the interval $]0,T]$. Then
$A^\e$, defined by \eqref{eq:opgen},  
satisfies \eqref{eq:commut} if and only if ${\tt V}$ is of the form
\eqref{def:V1}. 
\end{prop}
\begin{rema}
We do not assume that $\phi$ solves the eikonal equation
\eqref{eq:eikon}. However, we will see in the proof that it is
essentially necessary. 
\end{rema}
\begin{rema}
Since from
Von Neumann equation, Heisenberg observables always satisfy
\eqref{eq:commut}, the above proposition implies that such an
observable can be written under the form \eqref{eq:opgen}, for some
functions $f$ and $\phi$, if and only if the potential ${\tt V}$ is of
the form \eqref{def:V1}.  
\end{rema}
\begin{proof}
We now only have to prove the ``only if'' part. 
Computations yield
\begin{equation}\label{eq:bracket}
\begin{split}
\left[ i\e \d_t +\frac{1}{2}\e^2\Delta -{\tt V}(x), A_j^\e(t)\right]=
&
f'(t)\d_j \phi +f(t)\d_{jt}^2 \phi + f(t)\d_j {\tt V} \\
 +\e\Big( -f'(t)&\d_j +f(t)\nabla_x(\d_j \phi)\cdot \nabla_x
+\frac{1}{2} f(t)\Delta (\d_j \phi)\Big). 
\end{split}
\end{equation}
This bracket is zero if and only if the terms in $\e^0$ and $\e^1$ are
zero. The term in $\e$ is the sum of an operator of order one and of
an operator of order zero. It is zero if and only if both operators
are zero. The operator of order one is zero if and only if
\begin{equation*}
f(t)\d_{jj}^2 \phi = f'(t)\ ,\ \ \d_{jk}^2 \phi \equiv 0 \
\textrm{ if } j\not = k\, .
\end{equation*}
In particular, $\d_{jj}^2 \phi$ is a function of time only,
independent of $x$, and we have
$$\frac{1}{2} f(t)\Delta (\d_j \phi)\equiv 0  \, .$$
%%Since $\phi$ solves the eikonal equation \eqref{eq:eikon}, the term
%%$\e^0$ also writes
From the above computations, the first two terms in $\e^0$ also write
\begin{equation*}
f'(t)\d_j \phi +f(t)\d_{jt}^2 \phi = \sum_{k=1}^n f(t)\d_k \phi
\d^2_{jk} \phi +  
f(t)\d_{jt}^2 \phi= f(t)\d_j \left(\d_t\phi + \frac{1}{2}|\nabla_x
\phi|^2\right) \, .
\end{equation*}
Canceling the term in $\e^0$ in \eqref{eq:bracket} therefore yields,
since $f$ is never zero on $]0,T]$, 
\begin{equation}\label{eq:lbm}
\d_j \left(\d_t\phi + \frac{1}{2}|\nabla_x
\phi|^2+ {\tt V}(x)\right) =0\, .
\end{equation}
Differentiating the above equation with respect to $x_k$ and $x_\ell$,
all the terms with $\phi$ vanish, since we noticed that the
derivatives of order at least three of $\phi$ are zero. We deduce that
for any triplet $(j,k,\ell)$, $\d_{jk\ell}^3 {\tt V} \equiv 0,$
that is, ${\tt V}$ is of the form \eqref{def:V1}.

Notice that since \eqref{eq:lbm} holds for any $j\in \{1,\ldots,n\}$,
there exists a 
function $\Xi$ of time only such that
\begin{equation*}
\d_t\phi + \frac{1}{2}|\nabla_x
\phi|^2+ {\tt V}(x)=\Xi (t)\, .
\end{equation*}
This means that $\phi$ is almost a solution to the eikonal equation
\eqref{eq:eikon}. Replacing $\phi$ by
$\widetilde\phi(t,x):=\phi(t,x) -\int_0^t 
\Xi(s)ds$ does not affect \eqref{eq:opgen},  and $\widetilde \phi$
solves \eqref{eq:eikon}.
%%$$f'(t)\d_j \phi -f(t)\sum_{k=1}^n \d_k \phi \d^2_{jk}\phi = 
%%\d_j \phi \left( f'(t)- f(t) \d^2_{jj}\phi\right)=0\ ,$$
%%from the above computations.  The proposition stems from the following
%%lemma. 
%%\begin{lem}\label{lem:d3}
%%Let ${\tt V}\in C^\infty(\R^n;\R)$. Let $\phi\in
%%C^4(]0,T[\times\R^n;\R)$ solve \eqref{eq:eikapp} for some function
%%$\Xi$. If for any triplet $(j,k,\ell)$, we
%%have $\d_{jk\ell}^3 \phi \equiv 0,$
%%then for every triplet $(j,k,\ell)$, $\d_{jk\ell}^3 {\tt V} \equiv 0,$
%%that is, ${\tt V}$ is of the form \eqref{def:V1}. 
%%\end{lem}
%%\begin{proof}[Proof of Lemma~\ref{lem:d3}]
%%Differentiate \eqref{eq:eikapp} with respect to $x_j$, $x_k$ and
%%$x_\ell$.
%%, we have
%%\begin{equation*}
%%\d_t \phi +\frac{1}{2}\sum_{s=1}^n(\d_s \phi)^2 =-{\tt V}(t,x)\ .
%%\end{equation*}
%%Differentiate with respect to $x_j$, this yields
%%\begin{equation*}
%%\d^2_{tj} \phi +\sum_{s=1}^n\d_s \phi\d^2_{sj} \phi =-\d_j{\tt V}(t,x)\ .
%%\end{equation*}
%%Differentiate with respect to $x_k$,
%%\begin{equation*}
%%\d^3_{tjk} \phi +\sum_{s=1}^n\left(\d^2_{sk} \phi\d^2_{sj} \phi +\d_s
%%\phi\d^3_{sjk} \phi \right)
%%=-\d^2_{jk}{\tt V}(t,x)\ . 
%%\end{equation*} 
%%From our assumption, this is also
%%\begin{equation*}
%%\d^3_{tjk} \phi +\sum_{s=1}^n\d^2_{sk} \phi\d^2_{sj} \phi 
%%=-\d^2_{jk}{\tt V}(t,x)\ . 
%%\end{equation*} 
%%Differentiating with respect to $x_\ell$ and using our assumption once
%%again completes the proof of Lemma~\ref{lem:d3}. 
%%\end{proof}
%%This completes the proof of the proposition.
\end{proof}
\subsection{Heisenberg observables for general subquadratic potentials}
We now suppose that ${\tt V}={\tt V}(t,x)$  satisfies
 Assumption~\ref{hyp:gen}. 
 Define the Heisenberg observable
$${\tt A}^\e(t)={\tt U}^\e(t)\frac{x}{\e}{\tt U}^\e(-t),$$
where the group ${\tt U}^\e$ is defined by \eqref{def:ttU}. 
The latter is in general not a differential operator, but a
pseudo-differential operator (Egorov theorem, see
e.g. \cite{Robert}). We saw that if $\tt V$ satisfies
Assumption~\ref{hyp:pot} however, then it is explicit. The drawback of
this approach is that we cannot assess the action of this operator on
nonlinear terms in general. The operator ${\tt A}^\e$ satisfies two of the
three properties we use to study the nonlinear problem:  
\begin{lem}\label{lem:Heisenberg}
The operator ${\tt A}^\e(t)$ satisfies the following properties.
\begin{itemize}
\item The commutation,
$$\left[{\tt A}^\e(t), i\e \d_t +\frac{1}{2}\e^2\Delta -{\tt V}(t,x)
\right]=0.$$
\item The modified Sobolev inequality. If $v\in\Sigma$, then for 
$2\leq r\leq \frac{2n}{n-2}$, there exists $C_r$ such that, for
$|t|\leq \delta$, 
$$\|v\|_{L^r}  \leq \frac{C_r}{|
t|^{\delta(r)}}\|v\|_{L^2}^{1-\delta(r)}
\|{\tt A}^\e(t)v\|_{L^2}^{\delta(r)}.$$
\end{itemize}
\end{lem}
\begin{proof}
The first point stems from the definition of
${\tt A}^\e(t)$. For the second, let $g^\eps(t,x)={\tt U}^\eps(-t)v(x)$. We
know that for any $f \in L^2\cap L^1$, 
$$\|f\|_{L^2}=\|{\tt U}^\eps(t)f\|_{L^2},$$
and for $|t|\leq \delta$, from (\ref{eq:solfond}), 
$$\|{\tt U}^\eps(t)f\|_{L^\infty}\lesssim |\eps t|^{-n/2}\|f\|_{L^1}.$$
Interpolating these two estimates yields,
$$\|{\tt U}^\eps(t)f\|_{L^r}\lesssim |\eps t|^{-\delta(r)}\|f\|_{L^{r'}},$$
therefore,
$$\|{\tt U}^\eps(t)g^\eps(t)\|_{L^r}\lesssim
|\eps t|^{-\delta(r)}\|g^\eps(t)\|_{L^{r'}}.$$
Let $\lambda >0$, and write,
\begin{equation*}
\|g^\eps(t)\|_{L^{r'}}^{r'}  = \int_{|x|\leq \lambda}|g^\eps(t,x)|^{r'} dx
+\int_{|x|> \lambda}|g^\eps (t,x)|^{r'} dx.
\end{equation*}
Estimate the first term by H\"older's inequality, 
$$\int_{|x|\leq \lambda}|g^\eps(t,x)|^{r'} dx
\lesssim \lambda^{n/p'} \left(\int_{|x|\leq \lambda}|g^\eps(t,x)|^{r'p}
dx\right)^{1/p},$$
and choose $p =2/r' (\geq 1)$. Estimate the second term by the same
H\"older's inequality,  after inserting the factor $x$ as
follows, 
\begin{equation*}
\begin{split}
\int_{|x|> \lambda}|g^\eps (t,x)|^{r'} dx & = \int_{|x|>
\lambda}|x|^{-r'} 
|x|^{r'} |g^\eps (t,x)|^{r'} dx \\
& \leq \left(\int_{|x|>
\lambda}| x|^{-r'p'}dx\right)^{1/p'} \left(\int_{|x|>
\lambda} | xg^\eps (t,x)|^2dx\right)^{1/p}\\
& \lesssim \lambda^{n/p'-r'}\|xg^\eps (t,x)\|_{L^2}^{2/p}. 
\end{split}
\end{equation*}
In summary, we have the following estimate, for any $\lambda >0$, 
\begin{equation}\label{eq:lambda}
\|g^\eps(t)\|_{L^{r'}} \lesssim \lambda^{n/(p'r')}\|g^\eps(t)\|_{L^2}
+ \lambda^{n/(p'r')-1}\|xg^\eps (t,x)\|_{L^2}.
\end{equation}
Notice that $n/(p'r')= \delta(r)$, and equalize both terms of the
right hand side of (\ref{eq:lambda}), 
$$\lambda =\frac{\|xg^\eps (t,x)\|_{L^2}}{\|g^\eps(t)\|_{L^2}}.$$
This yields,
\begin{equation*}
\|g^\eps(t)\|_{L^{r'}} \lesssim \|g^\eps(t)\|_{L^2}^{1-\delta(r)}
\|xg^\eps (t,x)\|_{L^2}^{\delta(r)}.
\end{equation*}
Therefore,
\begin{equation*}
\begin{split}
\|{\tt U}^\eps(t)g^\eps(t)\|_{L^r} &\lesssim
|\eps t|^{-\delta(r)}\|g^\eps(t)\|_{L^2}^{1-\delta(r)}
\|xg^\eps (t,x)\|_{L^2}^{\delta(r)}\\
& \lesssim
|t|^{-\delta(r)}\|g^\eps(t)\|_{L^2}^{1-\delta(r)}
\left\|\frac{x}{\eps}g^\eps (t,x)\right \|_{L^2}^{\delta(r)}.
\end{split}
\end{equation*}
Back to $v$, this completes the proof of the lemma, since
${\tt U}^\eps(t)$ is unitary on $L^2$.
\end{proof}

\subsection{A partial result for general subquadratic potentials} 
To conclude, we prove that the asymptotics \eqref{eq:asym1} still
holds if ${\tt V}$ satisfies Assumption~\ref{hyp:gen}.
\begin{prop}
Let ${\tt V}$ satisfying Assumption~\ref{hyp:gen}, such that ${\tt V}$
is continuous at $(t,x)=(0,0)$, with ${\tt
V}(0,0)=0$. Suppose that
Assumption~\ref{hyp:nl} is 
satisfied. Then for any $\Lambda >0$, the following holds:\\
1. There exists $\e(\Lambda)>0$ such that for $0<\e\leq
   \e(\Lambda)$, the initial value problem
\begin{equation}\label{eq:pbgen}
\left\{
\begin{split}
i\e \d_t u^\e +\frac{1}{2}\e^2 \Delta u^\e &= {\tt V}(t,x)u^\e
+ \e^{n\si} |u^\e|^{2\si} u^\e\ ,\\
u^\e_{\mid t=0}&= \frac{1}{\e^{n/2}}\varphi\left(
\frac{x}{\e}\right) \ ,
\end{split}
\right.
\end{equation}
has a unique solution $u^\e \in
   C([-\Lambda\e,\Lambda\e];\Sigma)$. \\
2. This solution satisfies the following asymptotics,
\begin{equation}\label{eq:asym1gen}
\begin{split}
\limsup_{\e \to 0}\sup_{|t|\leq \Lambda \e}\Big(&
\left\|u^\e(t)-v^\e(t) \right\|_{L^2}+
\left\|\e\nabla_x u^\e(t)-\e\nabla_x v^\e(t) \right\|_{L^2}\\
& +
\left\|\frac{x}{\e}u^\e(t)- \frac{x}{\e}v^\e(t) \right\|_{L^2}\Big)
   =0\ ,
\end{split}
\end{equation}
where $v^\e$ is given by \eqref{eq:v}. 
\end{prop}
\begin{proof}
The proof mimics the approach used in Section~\ref{sec:cl}, except
that we do not use intermediary operators such as $A_{\ell,j}^\e$. 
Denote $w^\eps = u^\eps- v^\eps$. It solves
\begin{equation}\label{eq:w}\left\{
\begin{split}
i\eps \d_t w^\eps +\frac{1}{2}\Delta w^\eps &= {\tt V}(t,x)w^\eps +
{\tt V}(t,x)v^\eps + 
\eps^{n\sigma} \left(|u^\eps|^{2\sigma}u^\eps -
|v^\eps|^{2\sigma}v^\eps \right),\\ 
w^\eps_{\mid t=0}&=0.
\end{split}\right.
\end{equation}
Obviously,
\begin{equation}\label{eq:rq}
\left||u^\eps|^{2\sigma}u^\eps -
|v^\eps|^{2\sigma}v^\eps \right| \lesssim \left( |v^\eps|^{2\sigma} +
|w^\eps|^{2\sigma} \right)|w^\eps|.
\end{equation}
We know that there exists $C_0$ such that for any $t$,
$$\|\e \nabla_x  v^\eps(t)\|_{L^2}\leq C_0\ .$$
Since $w^\eps_{\mid t=0}=0$ and $w^\eps \in C(0,t^\eps;\Sigma)$ for
some $t^\eps >0$ (Proposition~\ref{prop:localex}), we have
\begin{equation}\label{eq:tantque1}
\| \e \nabla_x w^\eps(t)\|_{L^2}\leq C_0\ ,
\end{equation}
for $t$ in some interval $[0,t_1^\eps]$. So long as
\eqref{eq:tantque1} holds, we can get energy estimates from
\eqref{eq:w}, proceeding as in Section~\ref{sec:cl} and using the
Gagliardo-Nirenberg inequality
\begin{equation*}
\|f\|_{L^{\underline r}}\leq C \e^{-\delta({\underline
r})}\|f\|_{L^2}^{1- \delta({\underline r})} \|\e \nabla_x
f\|_{L^2}^{\delta({\underline r})} \ .
\end{equation*}
Notice that we have,
\begin{equation*}
\left[i\e \d_t +\frac{1}{2}\e^2\Delta -{\tt V}(t,x), \e
\nabla_x\right]= \e \nabla_x{\tt V}(t,x)\ ;\ \ \left[i\e \d_t
+\frac{1}{2}\e^2\Delta -{\tt V}(t,x), 
\frac{x}{\e}\right]=\e \nabla_x\ .
\end{equation*}
Proceeding as in Section~\ref{sec:cl} yields, 
\begin{equation}\label{eq:wL2}
\|w^\eps\|_{L^\infty(0,t;L^2)}\leq C(\Lambda)
\eps^{-1}\|{\tt V}(s,x)v^\eps\|_{L^1(0,t;L^2)}\ ,  
\end{equation}
along with
\begin{equation*}
\begin{split}
\|\eps \nabla_x w^\eps\|_{L^\infty(0,t;L^2)} \leq C(\Lambda)\Big( &
\|\nabla_x {\tt V}(s,x) w^\eps\|_{L^1(0,t;L^2)}
 + \left\|\nabla_x\left({\tt V}(s,x)v^\eps\right)
\right\|_{L^1(0,t;L^2)}\Big) ,\\
\left\|\frac{x}{\e} w^\eps\right\|_{L^\infty(0,t;L^2)} \leq C(\Lambda)
\Big( & \| \nabla_x w^\eps\|_{L^1(0,t;L^2)}+ 
\eps^{-2}\|x{\tt V}(s,x)v^\eps\|_{L^1(0,t;L^2)}\Big).
\end{split}
\end{equation*}
In particular, so long as
\eqref{eq:tantque1}
holds, with $|t|\leq \Lambda \eps$, 
\begin{equation*}
\begin{split}
\|w^\eps  \|_{L^\infty(0,t;L^2)} + \|\eps \nabla_x  
w^\eps & \|_{L^\infty(0,t;L^2)}  + \left\|\frac{x}{\e} 
w^\eps\right\|_{L^\infty(0,t;L^2)}
\leq\\
\leq  e^{C(\Lambda)\Lambda}\Bigl( &
\eps^{-1}\|{\tt V}(s,x)v^\eps\|_{L^1(0,t;L^2)}
 +
\left\|\nabla_x\left({\tt V}(s,x)v^\eps\right)\right\|_{L^1(0,t;L^2)}\\
&+ \eps^{-2}\|x{\tt V}(s,x)v^\eps\|_{L^1(0,t;L^2)}\Bigr).
\end{split}
\end{equation*}
Now,
\begin{equation*}
\begin{split}
\eps^{-1}\|{\tt V}(s,x)v^\eps\|_{L^1(0,t;L^2)} & =
\eps^{-1}\left\|{\tt V}(s,x)\frac{1}{\eps^{n/2}}\psi \left(\frac{s}{\eps},
\frac{x}{\eps}\right)\right \|_{L^1(0,t;L^2)} \\
& = \eps^{-1}\left\|{\tt V}(s,\eps x)\psi \left(\frac{s}{\eps},
x\right)\right \|_{L^1(0,t;L^2)}\\
& = \left\|{\tt V}(\eps s,\eps x)\psi (s,x)
\right \|_{L^1(0,t/\eps;L^2)}\\
& \leq \left\|{\tt V}(\eps s,\eps x)\psi (s,x)
\right \|_{L^1(0,\Lambda ;L^2)}.
\end{split}
\end{equation*}
Notice that for $|t|\leq \delta$, 
$$ |{\tt V}(t,x)|\lesssim 1+x^2.$$
From Lebesgue's dominated convergence theorem (${\tt V}$ is continuous at
the origin and ${\tt V}(0,0)=0$) and Lemma~\ref{lem:fur},
it follows, up to approximating $\varphi$ in ${\mathcal S}(\R^n)$ as
in Section~\ref{sec:cl}, 
$$\left\|{\tt V}(\eps s,\eps x)\psi (s,x)
\right \|_{L^1(0,\Lambda ;L^2)}\Tend \eps 0 0.$$
Similarly, 
$$\left\|\nabla_x\left({\tt V}(s,x)v^\eps\right)
\right\|_{L^1(0,\Lambda \eps;L^2)} 
+ \eps^{-2}\|x{\tt V}(s,x)v^\eps\|_{L^1(0,\Lambda \eps ;L^2)}\Tend \eps 0 0.$$
Therefore \eqref{eq:tantque1} remains valid up to time $t=\Lambda
\eps$, provided that $\eps$ is sufficiently small ($0< \eps \leq
\eps(\Lambda)$). This completes the proof of the proposition.
\end{proof}

\end{document}